\documentclass[11pt]{amsart}

\usepackage{amssymb,amsfonts,amsbsy}
\usepackage{amsmath,amsthm}
\usepackage[active]{srcltx}		
\usepackage{pdfsync}					
\newcommand{\bg}{\begin{equation}}
\newcommand{\ed}{\end{equation}}
\newcommand{\bga}{\begin{eqnarray}}
\newcommand{\eda}{\end{eqnarray}}
\newcommand{\pf}{\noindent\textbf{Proof.~}}

\newcommand{\Ff}{{\mathcal F}}

\newcommand{\Ss}{{\mathcal{S}}}

\newcommand{\Yy}{{\mathcal{Y}}}
\newcommand{\Xx}{{\mathcal{X}}}
\newcommand{\tXx}{\widetilde{\mathcal{X}}}
\newcommand{\tYy}{\widetilde{\mathcal{Y}}}

\def\endProof{{\hfill$\Box$}}
\newcommand\R{{\mathbb{R}}}
\renewcommand\P{{\mathbb{P}}}
\newcommand\N{{\mathbb{N}}}
\newcommand\C{{\mathbb{C}}}

\newcommand\K{{\mathbb{K}}}

\newtheorem{theorem}{Theorem}[section]
\newtheorem{proposition}[theorem]{Proposition}
\newtheorem{lemma}[theorem]{Lemma}
\newtheorem{corollary}[theorem]{Corollary}

\theoremstyle{definition}

\def\cbdu{\par{\raggedleft$\Box$\par}}

\theoremstyle{remark}
\newtheorem{remark}[theorem]{Remark}

\numberwithin{equation}{section}

\def\R{\mathbb{R}}
\def\C{\mathbb{C}}
\newcommand\Proof{\noindent{\bf Proof.~}}
\def\endProof{{\hfill$\Box$}}

\renewcommand{\P}{\mathbb{P}}

\begin{document}

\title[Decay and growth of solutions of Boussinesq (\today)]
{Large time decay and growth for solutions of a viscous Boussinesq system}

\author{Lorenzo Brandolese and Maria E. Schonbek}

\address{L. Brandolese: Universit\'e de Lyon~; Universit\'e Lyon 1~;
CNRS UMR 5208 Institut Camille Jordan,
43 bd. du 11 novembre,
Villeurbanne Cedex F-69622, France.}
\email{brandolese{@}math.univ-lyon1.fr}
\urladdr{http://math.univ-lyon1.fr/$\sim$brandolese}
\address{M.~E.~Schonbek: Department of Mathematics, UC Santa Cruz, Santa Cruz, CA 95064,USA}
\email{schonbek@math.ucsc.edu}

\thanks{The work of L. Brandolese and M. Schonbek were partially supported by FBF Grant SC-08-34.
The work of M. Schonbek was also partially supported by NSF Grants DMS-0900909 and 
grant FRG-09523-503114.}


\date{\today}


\begin{abstract} In this paper we analyze the decay  and the growth for large time of weak and strong
solutions to the three-dimensional viscous Boussinesq system. We show that generic solutions blow up as $t\to\infty$
in the sense that the energy  and  the $L^p$-norms of
the velocity field grow to infinity for large time for $1\le p<3$.
 In the case of strong solutions we provide sharp estimates  both from above and from below and explicit asymptotic profiles.
We also show that solutions arising from $(u_0,\theta_0)$ with zero-mean for the initial temperature
$\theta_0$ have a special behavior as $|x|$ or $t$ tends to infinity: contrarily to the generic case, their energy dissipates to zero for large time.
 \end{abstract}

\maketitle

\section{Introduction }

In this paper we address the problem of the heat transfer inside viscous incompressible
flows in  the whole space $\R^3$.
Accordingly with  the Boussinesq approximation, we neglect the variations of the density
in the continuity equation and the local heat source due to the viscous dissipation.
We rather take into account the variations of the temperature by putting an
additional vertical  buoyancy force term in the equation of the fluid motion.

This leads us to the Cauchy problem for the Boussinesq system
\begin{equation}
\label{B} 
\left\{
\begin{aligned}
 &\partial_t \theta +u \cdot \nabla \theta =  \kappa\Delta \theta\\
 &\partial_t u +u\cdot\nabla u+\nabla p= \nu\Delta u+\beta\theta e_3\\
  &\nabla\cdot u=0\\
 &u|_{t=0}=u_0,\;\,\theta|_{t=0}=\theta_0.
\end{aligned}
\right.
\qquad x\in \R^3, t\in \R_{+}
\end{equation}
Here $u\colon\R^3\times\R^+\to \R^3$ is the velocity field. The scalar fields 
$p\colon \R^3\times\R^+\to \R$ and $\theta\colon\R^3\times\R^+\to\R$ denote respectively 
the pressure and the  temperature of the fluid. Moreover, $e_3=(0,0,1)$, and $\beta\in \R$
is a physical constant.
For the decay questions that we address in this paper, it will be
important to have strictly positive
viscosities in both equations: $\nu,\kappa>0$.
By rescaling the unknowns, we can and do assume, without loss of generality, that $\nu=1$ and $\beta=1$.
To simplify the notation, from now on we take the thermal diffusion coefficient $\kappa>0$ such that $\kappa=1$.

Like for the Navier--Stokes equations, obtained as a particular case from~\eqref{B} putting
$\theta\equiv0$,  weak solutions to~\eqref{B} do exist, but their uniqueness is not known.
The global existence of weak solutions, or strong solutions in the case of small data has been 
studied by several authors. See, {\it e.g.\/} \cite{AbidiHmidi07}, \cite{dongochae},  \cite{DaP1}, \cite{DaP2}, \cite{SawT}.
Conditional regularity results for weak solutions (of Serrin type) can be found in~\cite{CdiB}.
The smoothness of solutions arising from large axisymmetric data is addressed in
\cite{AbidiHmidiKer} and \cite{HmidiRou10}.
Further regularity issues on the solutions have been discussed  also  \cite{FanZhou09, FefCordoba}.

\medskip
The goal of this paper is to study in which way the variations of the temperature affect the
asymptotic behavior of the velocity field.
We point out that several different models are known in the literature under the name
of ``viscous (or dissipative) Boussinesq system''.
The asymptotic behaviour of viscous Boussinesq systems of different
nature have been recently addressed, {\it e.g.\/}, in~\cite{AbouAG09, ChenG09}.
But the results therein cannot be compared with ours.

\medskip
Only few works are devoted to the study of the large  time behavior of solutions to~\eqref{B}.
See \cite{FerVil06, KarPri08}.
These two papers deal with self-similarity issues and stability results for solutions in critical spaces
(with respect to the scaling).
On the other hand, we will be mainly concerned with {\it instability\/} results for the energy norm,
or for other subcritical spaces, such as $L^p$, with $p<3$.

\medskip

A simple energy argument shows that weak solutions arising from data $\theta_0\in L^1\cap L^2$ and $u_0\in L^2_{\boldsymbol{\sigma}}$
satisfy the estimates
$$\|u(t)\|_2\le C(1+t)^{1/4}$$ 
and 
$$\|\theta(t)\|_2\le C(1+t)^{-3/4}.$$
The above estimate for the temperature looks optimal, since the decay agrees with that of the heat kernel.
On the other hand the optimality of the estimate for the velocity field is not so clear.

For example, in the particular case $\theta_0=0$, the system boils down to the  Navier--Stokes equations
and in this simpler case one can improve the bound for the velocity  into $\|u(t)\|_2\le \|u_0\|_2$. In fact, $\|u(t)\|_2\to0$ for large time
by a result of Masuda~\cite{Mas84}.
Moreover, in the case of Navier--Stokes the decay of $\|u(t)\|_2$ agrees with the 
$L^2$-decay of the solution of the heat equation.
See \cite{Schon84, KM, Wie87} for a more precise statement.

\medskip
The goal of the paper will be to show that the estimate of weak solutions $\|u(t)\|_2\le C(1+t)^{1/4}$
can be improved if and only if the initial temperature has zero mean.
To achieve this, we will  establish the validity of the  corresponding lower bounds
for a class of {\it strong solutions\/}.

In particular, this means that
very  nice data (say, data that are smooth, fast decaying and ``small'' in some strong norm)
give rise to solutions that become large as $t\to\infty$:
our results imply the growth of the energy
 for strong solutions~:
\begin{equation} 
 \label{lbu2}
  c(1+t)^{1/4}\le \|u(t)\|_2\le C(1+t)^{1/4},
\qquad t>\!\!\!>1.
\end{equation}
The validity of the lower bound in~\eqref{lbu2}  (namely, the condition $c>0$), will be ensured
whenever the initial temperature is sufficiently decaying but
$$\int \theta_0\not=0.$$

\medskip
We feel that is important to point out here an erratum to the paper~\cite{GuoGuang94}.
Unfortunately, the lower bound in~\eqref{lbu2} contredicts a result in~\cite[Theorem~2.3]{GuoGuang94},
where the authors claimed that $\|u(t)\|_2\to 0$ under too general assumptions, weaker than those leading to our
growth estimate.
The proof of their theorem (in particular, of inequalities (5.6) and (5.8) in \cite{GuoGuang94}) can be fixed 
by putting different conditions on the data, including  $\int\theta_0=0$.
This is essentially what we will do in part (b) of our Theorem~\ref{th:belo} below.
Similarity, the statement of Theorem~2.4 in \cite{GuoGuang94} contredicts our lower bound~\eqref{sda2} below
(inequalities (5.18)--(5.20) in their proof do not look correct).
This will be also corrected by our Theorem~\ref{th:belo}.
We would like  to give credit to the paper~\cite{GuoGuang94} 
(despite the above mentioned errata), because we got from there
inspiration for our results of Sections~\ref{sec:three} and~\ref{secsix}.

\medskip
Our main tool for establishing the the lower bound will be the derivation of exact pointwise asymptotic profiles of solutions
in the parabolic region $|x|>\!\!\!>\sqrt t$.
This will require a careful choice of several function spaces in order to
obtain as much information as possible on the pointwise behavior of the
velocity and the temperature.
A similar method has been applied before by the first author in
\cite{BraV07} in the case of the Navier--Stokes equations,
although the relevant estimates were performed there in a different functional setting.

Even though several other methods developped for  Navier--Stokes
could be effective for obtaining estimates from below
(see, e.g.,~\cite{choe-jin, GW2, MiyS}), our analysis has the advantage of putting in evidence
some features that are specific of the Boussinesq system: in particular,  the different behavior
of the flow when $|x_3|\to\infty$ or when $\sqrt{x_1^2+x_2^2}\to\infty$, due to 
the verticality of the bouyancy forcing term $\theta e_3$ (see Theorem~\ref{theor1} below).
Moreover, the analysis of solutions in the region $|x|>\!\!\!>\sqrt t$ and our use of weighted spaces
completely explains the phenomenon of the energy growth:
the variations of the temperature push the fluid particles in the far field; even though in any bounded region
the fluid particles slow down as $t\to\infty$ (this effect is measured {\it e.g.\/} by the decay of the $L^\infty$-norm established
in Proposition~\ref{prop2}),  large portions of fluid globally carry an increasing energy during the evolution.
Our result thus illustrates the physical limitations of the Boussinesq approximation, at least  for the study of heat convection inside fluids 
filling  domains where Poincar\'e's inequality is not available, such as the whole space.

\medskip
In fact, our method applies also to weighted $L^p$-spaces, so let us introduce the 
weighted norm
$$\|f\|_{L^p_r} = \biggl(\int |f(x)|^p(1+|x|)^{pr}\,dx\biggr)^{1/p}.$$
Then we will show that  strong solutions starting from suitably small and well decaying data
satisfy, for
$t>0$ large enough,
\begin{equation}
\label{sda}
c(1+t)^{\frac{1}{2}(r+\frac{3}{p}-1)}\le \|u(t)\|_{L^p_r}
\le C(1+t)^{\frac{1}{2}(r+\frac{3}{p}-1)},
\end{equation}
for all 
\begin{equation*}
\label{para} r\ge 0,\quad 1<p<\infty,\qquad r+\frac{3}{p}<3.
\end{equation*}
As before in~\eqref{lbu2}, these lower bounds hold true with a constant $c>0$ as soon as the initial temperature has
non-zero mean.
Notice that in this case the $L^p$-norms asymptotically blow up for large time if and only if $p<3$.

The above restriction $r+\frac{3}{p}<3$ on the parameters is optimal.
Indeed, under the same conditions yielding to~\eqref{sda} we will also prove that
when
$ r\ge 0$, $1\le p<\infty$ and  $r+\frac{3}{p}\ge 3$, then one has
$\|u(t)\|_{L^p_r}=\infty$, for all $t>0$.
The fact that the $L^p_r$-norm becomes infinite instanteneously  for this 
range of the parameters is related to  that  the velocity field immediately spatially spreads out
 and cannot decay faster than $|x|^{-3}$ as $|x|\to\infty$ for $t>0$,
 and this even if $u_0\in C^\infty_0(\R^3)$.

\bigskip
As we observed before, the lower bound in~\eqref{sda} brakes down when 
$\int\theta_0=0$.
In such case, our decay estimates can be improved.
We will establish the  upper bound for  {\it weak solutions\/}
\begin{equation}
\label{ades}
\|u(t)\|_2\le C'(1+t)^{-1/4}
\end{equation}
by the Fourier splitting method. This method was first introduced in~\cite{Schon84}.
We will need no smallness assumption on~$u_0$ to prove~\eqref{ades}.
We have to put however  a smallness condition of the form
\begin{equation}
\label{getrid}
\|\theta_0\|_1<\varepsilon_0.
\end{equation}
We do not know if it is possible to get rid  of~\eqref{getrid} to establish~\eqref{ades}.
Such smallness condition however looks natural as it respects the natural scaling invariance of the system~\eqref{B}.

\medskip
As before, the decay estimate~\eqref{ades} is optimal for generic solutions satisfying 
$\int\theta_0=0$. 
Indeed, when we start from  localized and small velocity, we can establish the upper-lower bounds
for strong solutions
\begin{equation}
\label{sda2}
c'(1+t)^{\frac{1}{2}(r+\frac{3}{p}-2)}\le \|u(t)\|_{L^p_r}
\le C'(1+t)^{\frac{1}{2}(r+\frac{3}{p}-2)},
\end{equation}
for all 
\begin{equation*}
\label{para-1} r\ge 0,\quad 1<p<\infty,\qquad r+\frac{3}{p}<4.
\end{equation*}
Similarily as before, the validity of the lower bound (the condition $c'>0$) now requires 
a non vanishing condition in the first moments of $\theta$.

\subsection{Notations}
\label{sec:notation}
We denote by $C_0^\infty$  the space of smooth functions with compact support.  
The $L^p$-will be denoted by $\|\cdot\|_p$ and in the case $p=2$ we will simply write $\|\cdot\|$. 
Moreover,
 \begin{equation*}
\mathcal{V}=\{ \phi \in C_0^\infty \,|\, \nabla\cdot \phi=0\},
\end{equation*}
$L^p_{\boldsymbol{\sigma}}$ denotes the completion of $\mathcal{V}$ under the norm $\|\cdot\|_p$,
and 
$V$ the closure of $\mathcal{V}$ in $H^1_0$

We will adopt the following convention for the Fourier transform of integrable functions:
$\Ff f(\xi)=\widehat{f}(\xi) = \int f(x)e^{-2\pi i x\cdot\xi}\, dx$.
Here and throughout the paper all integral without integration limits are over the whole $\R^3$.

The notations $L^{p,q}$ and $L^p_r$ have a different meaning.
For $1<p<\infty$ and $1\le q \le \infty$, $L^{p,q}$ denotes the classical Lorentz space.
For the definition and the basic inequalities concerning Lorentz spaces (namely the generalisation
of the straightforward $L^p$-$L^q$ H\"older and Young convolution inequalities) the reader can refer to
\cite[Chapter 2]{Lem02}.
Here we just recall that $L^{p,p}$ agrees with the usual Lebesgue space $L^p$, and
$L^{p,\infty}$ agrees with the weak Lebesgue space $L^{p}_w$
(or Marcinkiewicz space)
 \[ L^{p}_w = \{ f:\R^n \to \C, \;\mbox{measurable}\; ,\, \|f\|_{L^p_w} <\infty\}.\]
The quasi-norm  
\begin{equation*}
 \|f\|_{L^p_w } = \sup_{t>0} t[ \lambda_f(t)]^{\frac{1}{p}} 
\end{equation*}
is equivalent to the natural norm on $L^{p,\infty}$, for $1<p<\infty$. Here as it is usual we defined 
\[ \lambda_f(s) = \lambda\{ x: f(x) >s\} \]
where $\lambda$ denotes the Lebesgue measure.
On the other hand, for $1\le p\le\infty$ and $0\le r\le\infty$, $L^p_r$ is the weighted Lebesgue space, consisting of the functions $f$
such that $(1+|x|)^q f\in L^p$.
Notice that the bold subscript $\boldsymbol{\sigma}$ introduced above is not a real parameter:
in the notation $L^p_{\boldsymbol{\sigma}}$, this subscript simply stands for ``solenoidal''.

 To denote general constants we use $C$ which may change from line to line.  In certain cases we will write $C(\alpha)$ to
 emphasize the constants dependence on $\alpha$.

We denote by  $e^{t\Delta}$ the heat semigroup. Thus, $e^{t\Delta}u_0=\int g_t(x-y)u_0(y)\,dy$, where
$g_t(x)=(4\pi)^{-3/2}e^{-|x|^2/(4t)}$ is the heat kernel.

We denote by $E(x)$ the fundamental solution of $-\Delta$ in $\R^3$.
The partial derivatives of $E$ are denoted by $E_{x_j}$, $E_{x_j,x_k}$, etc.

\subsection{Organization of the paper}
All the main results are  stated without proof in Section~\ref{sec:two}.
The statement of our theorems are splitted into two parts: part (a) is devoted to the properties
of solutions in the  general case (where the integral $\int \theta_0$ is not necessarily zero).
In Part (b) of our theorems we are concerned with the special case $\int\theta_0=0$.

The rest of the paper is organized as follows. 
Sections~\ref{sec:three}-\ref{secsix} are
devoted to the proof of our results on weak solutions and
Sections~\ref{sec:seven}-\ref{sec:eight} to strong solutions.
In Section~\ref{sec:ee8} we collect a few technical remarks.

\section{Statement of the main results}
\label{sec:two}

\subsection{Results on weak solutions}
\label{sec:two.one}

 We will consider  weak solutions to the viscous Boussinesq equations and establish existence
 and natural decay estimates 
in $L^p$ spaces,  $1\leq p< \infty$, for the temperature, together with bounds for the growth of the velocities.
Such estimates rely on the fact that in both equations of the system~\eqref{B} we have a diffusion term. They complete
those of obtained in \cite{DP} where there was no temperature diffusion.

We start recalling the basic existence result of weak solutions to the system~\eqref{B}.
See, {\it e.g.\/}, \cite{CdiB, DP}.

\begin{proposition}
 \label{th:existen}
Let $(\theta_0,u_0) \in  L^2 \times L^2_{\boldsymbol{\sigma}} $. There exists a weak solution $(\theta, u)$ of the Boussinesq system (\ref{B}), continuous from $\R^+$ to $L^2$ with the weak topology,
 with  data $(u_0,\theta_0) $ such that, for any $T>0$,
$$\theta \in L^2(0,T; H^1) \cap L^{\infty}(0,T; L_{\boldsymbol{\sigma}}^2),\qquad
  u \in L^2(0,T; V) \cap L^{\infty}(0,T; L_{\boldsymbol{\sigma}}^2).$$
Such solution satisfies, for all $t\in[0,T]$, the energy inequalities
\begin{equation}
 \label{eit} 
\|\theta(t)\|^2 
+2\int_0^t \|\nabla \theta(s)\|^2\,ds \leq  \|\theta_0\|^2.
\end{equation}
and
\begin{equation}
\label{eiu}
\|u(t)\|^2+2\int_0^t\|\nabla u(s)\|^2\,ds
  \le +C\Bigl( \|u_0\|+t^2\|\theta_0\|^2\Bigr)
\end{equation}
for all $t\ge0$, and some absolute constant $C>0$.
\end{proposition}

\medskip
One can improve the growth estimate on the velocity
as soon as $\theta_0$ belongs to some $L^p$ space, with $p<2$.
For simplicity we will
consider only the case $\theta_0\in L^1\cap L^2$.

Moreover, it is natural to ask under which supplementary conditions on the initial data one can 
insure that the energy of the fluid $\|u(t)\|^2$ remains uniformly bounded.
Theorem~\ref{th:belo} provides an answer. 
Beside a  smallness assumption on $\|\theta_0\|_1$, we need to assume that $\int\theta_0=0$.

In addition, it is possible to prove that $\|u(t)\|$ not only remains uniformly bounded,
but actually  decays at infinity (without any rate).
Explicit decay rates for $\|u(t)\|$ can also be prescribed provided 
the linear part $e^{\Delta} u_0$ decays at the appropriate rate.

\begin{theorem}
 \label{th:belo}
\begin{enumerate}
 \item[(a)]
Let $(\theta_0,u_0) \in  L^2 \times L^2_{\boldsymbol{\sigma}} $.
Under the additional condition $\theta_0\in L^1$ the estimates on the weak solution
constructed in Proposition~\ref{th:existen} can be improved into
\begin{equation}
\label{CCC1}
\begin{split}
& \|\theta(t)\|^2 \le C(t+1)^{-\frac{3}{2}},\\
& \|u(t)\|^2   \leq C(t+1)^{\frac{1}{2}}, \\
\end{split}
\end{equation}
Moreover, if  $\theta_0 \in   L^1\cap L^p,\; \mbox{for some} \;1\leq p< \infty $, then 
$$\|\theta(t)\|_p \le  C(p)(t+1)^{-\frac{3}{2}(1-\frac{1}{p})}.$$

\item[(b)]
(The $\int\theta_0=0$ case)
In this part we additionally assume $\theta_0\in L^1_1$ and $\int\theta_0=0$.
Then there exists an absolute constant $\varepsilon_0>0$ such that if
\begin{equation}
 \label{varesss}
\|\theta_0\|_1<\varepsilon_0
\end{equation}
then the  weak solution of  the Boussinesq system~\eqref{B} constructed in Proposition~\ref{th:existen}
satisfies, for some constant $C>0$ and all $t\in\R^+$,
 \begin{align}
\label{bestt}
\|\theta(t)\| ^2\le& C(1+t)^{-\frac{5}{2}}
\end{align}
and
\begin{align}
\|u(t)\|^2 \to0 \qquad \hbox{as $t\to\infty$}.
\end{align}
Moreover under the additional condition~$u_0\in L^{3/2}\cap L^2_{\boldsymbol{\sigma}}$,  we have
\begin{align}
\label{bestu}
\|u(t)\|^2 \leq &C(1+t)^{-1/2}.
\end{align}
\end{enumerate}
\end{theorem}

The proof of part (a) of Theorem~\ref{th:belo} is straightforward. Part (b) is more subtle:
its proof relies on an inductive argument: assuming that the approximate velocity $u^{n-1}$ grows in $L^2$ at most like
$t^{1/8}$ (which is actually better than provided by Proposition~\ref{th:existen}) 
we can prove that the same growth estimate  remains valid for $u^n$.
A smallness condition on $\theta_0$ is needed to insure that in the inequality $\|u^n(t)\|\le C(1+t)^{1/8}$ the constant 
can by taken independently on~$n$.
With this improved control on the growth of the velocity, applying  the Fourier splitting method
(introduced in~\cite{Schon84}) and a few  boot-strapping 
we can improve our estimates up to the rates given in~\eqref{bestt} and \eqref{bestu}.

\begin{remark}
Solutions of the Boussinesq system~\eqref{B} with energy decaying faster
than~$t^{-1/2}$ might exist, but they are likely to be highly non-generic.
Indeed, it seems difficult to construct such solutions  assuming only that the data
belong to suitable function spaces (possibly with small norms).
The main obstruction is that one would need stringent cancellations properties
on the data, which however turn out to be  non-invariant under the Boussinesq flow.
A possible way to obtain such fast decaying  solutions would be to start with data satisfying
some special rotational symmetries, as those described in~\cite{Bra04iii}.
\end{remark}

\subsection{Results on strong solutions}
\label{sec:two.two}

The best way to prove the optimality of the estimates contained in Theorem~\ref{th:belo}
is to establish the corresponding lower bound estimates at least for a subclass of solutions.
For the study of the estimates from below
we will limit our considerations to a class of strong solutions.
This is not a real restriction as  lower bound estimates established 
for  solutions emanating from well localized, smooth and small data
are expected to remain valid in the larger class of weak solutions.
Studying strong solutions has also the advantage of better putting in evidence some interesting properties
specific of the Boussinesq system, such as the influence of the vertical 
 buoyancy force on the pointwise behavior of the fluid in the far-field.

\medskip
The existence of strong solutions to the system~\eqref{B} will be insured by 
a fixed point theorem in  function spaces invariant under 
the natural scaling of the equation. 
Thus, if $u\in \Xx$ , where $\Xx$  is a Banach space
to be determined, we want to have, for all $\lambda>0$,
\[\|u_\lambda\|_\Xx=\|u\|_\Xx, \qquad\hbox{where}\qquad
u_\lambda(x,t)=\lambda u(\lambda x,\lambda^2t)
\]
 is the rescaled velocity.
 A suitable  choice for norm of the space~$\Xx$, inspired by~\cite{CazDW}, is
\begin{equation}
\label{n21}\|u\|_{\mathcal{X}}=\,\,\,\,\!
\hbox{ess\,\,\,\,\,\,\,}
\!\!\!\!\!\!\!\!\!\sup_{x\in\R^3,\, t>0}\,\,\,\, (\sqrt t+|x|)\bigl|u(x,t)\bigr|.
\end{equation}
This choice for~$\Xx$ is quite natural. Indeed, whenever $|u_0(x)|\le C|x|^{-1}$, the linear evolution $e^{t\Delta}u_0$ belongs to~$\Xx$ and this paves the way for the  application of the fixed point
theorem in such space.

More precisely,  we define $\mathcal{X}$ as the Banach space of all locally integrable divergence-free vector fields $u$ such that
$\|u\|_{\mathcal{X}}<\infty$, and continuous with respect to $t$ in the following usual sense:
$u(t)\to u(0)$ in the distributional sense as $t\to0$ and $\hbox{ess}\,\sup_{x\in\R^3}|x|\,|u(x,t)-u(x,t')|\to0$ as $t\to t'$
if $t'>0$.

In the same way, if $\theta$  belongs to a  Banach space $\Yy$, we want to have
 $$\|\theta_\lambda\|_\Yy=\|\theta\|_\Yy, \qquad\hbox{where}\qquad
\theta_\lambda(x,t)=\lambda^3 \theta(\lambda x,\lambda^2t)$$
is the rescaled temperature.
 We then define a Banach space  $\mathcal{Y}$ of scalar functions through
the norm
\begin{equation}
\|\theta\|_{\mathcal{Y}}=
\|\theta\|_{L^\infty_t(L^1)}+
\hbox{ess\,\,\,\,\,\,\,}
\!\!\!\!\!\!\!\!\!\sup_{x\in\R^3,\, t>0}\,\,\,\, (\sqrt t+|x|)^{3}\bigl|\theta(x,t)\bigr|
\end{equation}
and the natural continuity condition on the time variable as before.

The starting point of our analysis will be the following proposition providing
a simple construction of  mild solutions $(u,\theta)\in \Xx\times\Yy$.
 We will refer to them also as {\it strong solutions\/}. Indeed,
one could prove that such solutions turn out to be smooth, as one  could check by adapting to
the system~\eqref{B} classical regularity criteria for the Navier--Stokes equations
like that of Serrin \cite{Ser62}. See the paper~\cite{CdiB}.
The smoothness of these solutions, however, plays no special role
in our arguments.

\begin{proposition}
 \label{prop1}
There exists an absolute constant $\epsilon>0$ such that if
\begin{equation}
\label{smalla}
 \|\theta_0\|_1<\epsilon, \qquad
\hbox{\rm ess}\,\sup_{x\in\R^3}|x|^3|\theta_0(x)|<\epsilon,
\qquad\hbox{\rm ess}\,\sup_{x\in\R^3}|x|\,|u_0(x)|<\epsilon
\end{equation}
where $u_0$ is a divergence-free vector field, then
there is a constant $C>0$ and a (mild) solution $(u,\theta)\in \mathcal{X}\times\mathcal{Y}$
of~\eqref{B}, such that
\begin{equation}
 \label{estsol}
\|u\|_{\mathcal{X}}\le C\epsilon \qquad\hbox{and}\qquad \|\theta\|_{\mathcal{Y}}\le C\epsilon.
\end{equation}
Moreover, these conditions define $u$ and $\theta$ uniquely.
\end{proposition}

Next Proposition shows it is possible to obtain better space-time decay estimates provided one starts with
suitably decaying data.

\medskip
\begin{proposition}
\label{prop2}
\begin{enumerate}
\item[(a)]
Let $u_0$ and $\theta_0$ as in Proposition~\ref{prop1}, and satisfying
the additional decay estimates, for some $1\le a<3$,   $b\ge3$, and a constant $C>0$,
\begin{equation}
 \label{adde}
\begin{split}
 |u_0(x)|\le C(1+|x|)^{-a},\\
 |\theta_0(x)|\le C(1+|x|)^{-b}.
\end{split}
\end{equation}
Then the solution constructed in Proposition~\ref{prop1} satisfies, for another constant $C>0$
independent on $x$ and $t$,
\begin{equation}
\label{deca1}
|u(x,t)|\le C\inf_{0\le \eta\le a} |x|^{-\eta}(1+t)^{(\eta-1)/2}
\end{equation}
and
\begin{equation}
\label{deca11}
|\theta(x,t)|\le C\inf_{0\le \eta\le b} |x|^{-\eta}(1+t)^{(\eta-3)/2}.
\end{equation}

%
%
\item[(b)]
(the  $\int \theta_0=0$ case)
Assume now $2\le a<4$, $a\not=3$, and $b\ge 4$ and let $u_0$ and $\theta_0$ satisfying the previous assumptions.
If, in addition,
\begin{equation}
 \label{via}
\int \theta_0=0 \qquad \hbox{and}\qquad \theta_0\in L^1_1,
\end{equation}
then the decay of $u$ and $\theta$ is improved as follows:
\begin{equation}
\label{deca2}
 \begin{split}
  |u(x,t)|\le C\inf_{0\le\eta\le a}(1+|x|)^{-\eta}(1+t)^{(\eta-2)/2},\\
 |\theta(x,t)|\le C\inf_{0\le \eta\le b}(1+|x|)^{-\eta}(1+t)^{(\eta-4)/2}.
 \end{split}
\end{equation}
\end{enumerate}
\end{proposition}

Recall that the fundamental solution of $-\Delta$ in~$\R^3$ is $E(x)=c|x|^{-1}$.
Thus, in the following asymptotic expansions, $\nabla E_{x_3}$ and 
$\nabla E_{x_h,x_3}$ are vectors
whose components are homogeneous functions of degree $-3$ and $-4$ respectively.
In particular $|\nabla E_{x_3}(x)|\le C|x|^{-3}$ and
$|\nabla E_{x_j,x_3}(x)|\le C|x|^{-4}$.

\medskip

We are now in the position of stating our main results on strong solutions.
The first theorem describes the asymptotic profiles of solutions in the parabolic region
$|x|>\!\!\!>\sqrt t$.
Roughly, it states that all sufficiently decaying solutions $(u,\theta)$ of~\eqref{B}
behave in such region like a potential flow.

\begin{theorem}\hfil
\label{theor1}
\begin{enumerate}
\item[(a)]
Let $a>\frac32$ and  $b>3$.
Let $(u,\theta)$ be a (mild) solution of~\eqref{B} satisfying the decay estimates~\eqref{deca1}-\eqref{deca11}.
Then the following profile for 
$u$ holds:
\begin{equation}  
\label{coroc}
u(x,t)=e^{t\Delta}u_{0}(x) +
\biggl(\int\theta_0\biggr)\,t\,\bigl(\nabla E_{x_3}\bigr)(x)\,+\,
\mathcal{R}(x,t)
\end{equation}
where  ${\mathcal{R}}(x,t)$
is a lower order term with respect to $\,t\,\nabla E_{x_3}(x)\,$ for $|x|>\!\!\!>\sqrt t$, namely,
\begin{equation}
 \lim_{\frac{|x|}{\sqrt t}\to\infty} \frac{{\mathcal{R}}(x,t)}{t |x|^{-3}}=0.
\end{equation}

\item[(b)] (the  $\int \theta_0=0$ case)
Assume now $a >2$ and  $b>4$. Assume also
that $\int \theta_0=0$.
Let $(u,\theta)$ be a solution satisfying the decay condition~\eqref{deca2}.
Then the following profiles for 
$u_j$  $(j=1,2,3)$ hold:
\begin{equation}
 \label{putz}
u_j(x,t)=e^{t\Delta}u_0(x)-\nabla E_{x_jx_3}(x)\cdot\biggl(\int_0^t\!\!\int y\,\theta(y,s)\,dy\,ds\biggr)+\widetilde{\mathcal{R}}(x,t)
\end{equation}
where  $\widetilde{\mathcal{R}}$
is a lower order term for $|x|>\!\!\!>\sqrt t>\!\!\!>1$, namely
\begin{equation}
 \lim_{t,\,\frac{|x|}{\sqrt t}\to\infty} \frac{\widetilde{\mathcal{R}}(x,t)}{t |x|^{-4}}=0.
\end{equation}
\end{enumerate}
\end{theorem}

The following remark should give a better understanding of the theorem.

\begin{remark}
 \label{rem-ullar}
\begin{enumerate}
 \item[(a)]
(The case $\int\theta_0\not=0$).
We deduce from the asymptotic profile~\eqref{coroc} the following:
when $|e^{t\Delta}u_0(x)|<\!\!\!<t|x|^{-3}$
(this happens, e.g., when we assume also
$|u_0(x)|\le C|x|^{-3}$ and $|x|>\!\!\!>\sqrt t>\!\!\!>1$)
and $\int\theta_0\not=0$ then
\begin{equation}  
\label{ullar}
u(x,t)\simeq \biggl(\int\theta_0\biggr)\,t\,\bigl(\nabla E_{x_3}\bigr)(x),
\qquad\hbox{for  $|x|>\!\!\!>\sqrt t$.}
\end{equation}
(The exact meaning of our notation and of statements~\eqref{ullar} and~\eqref{ullar2}
below  is made precise in the proof).
\item[(b)]
(The $\int\theta_0=0$ case)
We deduce from the profile~\eqref{putz} the following:
when
$|e^{t\Delta}u_0(x)|<\!\!\!<t|x|^{-4}$
(this happens, e.g., when we assume also $|u_0(x)|\le C|x|^{-4}$ and $|x|>\!\!\!>\sqrt t>\!\!\!>1$)
then
\begin{equation}
\label{ullar2}
u_j(x,t)\simeq 
-\nabla E_{x_jx_3}(x)\cdot\biggl(\int_0^t\!\!\int y\,\theta(y,s)\,dy\,ds\biggr)
\qquad
\hbox{for  $|x|>\!\!\!>\sqrt t>\!\!\!>1$.}
\end{equation}
\end{enumerate}
\end{remark}

\medskip
A remarkable consequence of the previous theorem is the following.

\begin{corollary}\hfil
\label{coro1}
\begin{enumerate}
\item[(a)]
Let $a>\frac32$, $b>3$ and let $(u,\theta)$ be a solution as in Part (a) of Theorem~\ref{theor1}.
Then for all $r,p$ such that
$$r\ge 0, \qquad 1<p<\infty,\qquad r+\frac{3}{p}<\min\{a,3\},$$
there exists $t_0>0$ such that the solution satisfies the upper and lower estimates
in the weighted-$L^p$-norm
\begin{equation}
\label{sdab}
\phi(|m_0|)\,\bigl(1+t\bigr)^{\frac{1}{2}(r+\frac{3}{p}-1)}\le \|u(t)\|_{L^p_r}
\le C'\bigl(1+t\bigr)^{\frac{1}{2}(r+\frac{3}{p}-1)}
\end{equation}
for all $t\ge t_0$.
Here, $m_0=\textstyle\int\theta_0$ and 
$\phi\colon\R^+\to\R$ is some continuous function such that
$\phi(0)=0$ and  $\phi(\sigma)>0$ if $\sigma>0$.
\item[(b)]
(the $\int\theta_0=0$ case)
Under the assumptions of the previous item,
with the stronger conditions
$a>2$,  $b>4$ and the additional zero mean condition $m_0=0$, let us set
$\widetilde{\boldsymbol{m}}=\liminf_{t\to\infty}\frac{1}{t}\left|\int_0^t\!\int y\theta(y,s)\,dy\,ds\right|$.
Then, for all $r,p$ such that
$$r\ge 0, \qquad 1<p<\infty,\qquad r+\frac{3}{p}<\min\{a,4\},$$
we have
\begin{equation}
\label{sdabz2}
\phi\bigl(\widetilde{\boldsymbol{m}}\bigr)\,\bigl(1+t\bigr)^{\frac{1}{2}(r+\frac{3}{p}-2)}\le \|u(t)\|_{L^p_r}
\le C'\bigl(1+t\bigr)^{\frac{1}{2}(r+\frac{3}{p}-2)}
\end{equation}
for another suitable continuous function $\phi\colon\R^+\to\R$ such that $\phi(0)=0$ and $\phi(\sigma)>0$ for $\sigma>0$.
\end{enumerate}
\end{corollary}

\begin{remark}
When $\int\theta_0\not=0$, we thus get by~\eqref{sdab}
the sharp large time behavior
$\|u(t)\|_{L^p_r}\simeq t^{\frac12(r+\frac3p-1)}$.

When $\int\theta_0=0$ and
$\widetilde{\boldsymbol{m}}\not=0$
we have the faster sharp decay
$\|u(t)\|_{L^p_r}\simeq t^{\frac12(r+\frac3p-2)}$.
The condition~$\widetilde{\boldsymbol{m}}\not=0$ is satisfied for generic solutions.
It prevents $\theta$ to have
oscillations at large times.
\end{remark}

\section{The mollified Boussinesq system and existence of weak solutions}
\label{sec:three}

The existence of weak solutions to the Boussinesq system
is well known, see \cite{CdiB}.
Their uniqueness, however, is an open problem.
Moreover, we do not know if {\it any\/} weak solutions satisfy the energy inequality and the decay estimates stated in
Proposition~\ref{th:existen}.
For this reason, we now briefly outline another construction of weak solutions,
which is well suited for obtaining all our estimates.

\medskip
We begin by introducing a mollified  Boussinesq system.
As the construction below is a straightworward adaptation of that of  Caffarelli, Kohn and Nirenberg,  \cite{CKN}, we will be rather sketchy. 
For completeness we recall the  definition  of the ``retarded mollifier" as given in~\cite{CKN}.
Let $\psi(x,t) \in C^{\infty}$ such that
\[
 \psi \geq 0,\;\; \int_0^{\infty}\!\!\! \int \psi dx\, dt =1,\;\; \mbox{supp}\,\psi \; \subset \{(x,t):|x|^2<t,1<t<2\}.
\]
For $T>0$ and $u \in L^2(0,T; L^2_{\boldsymbol{\sigma}})$, let  $\tilde{u}:\R^3\times \R \to \R^3$ be 
\[ 
\tilde{u} = 
\begin{cases}
u(x,t) &\mbox{if}\; (x,t) \in \R^3\times (0,T),\\
0&\mbox{otherwise.}
\end{cases}
\]
Let $\delta=T/n$. We set
\[ \Psi_{\delta}(u)(x,t) = \delta^{-4}\int_{\R^4}
 \psi\left(\frac{y}{\delta},\frac{\tau}{\delta}\right)\tilde{u}(x-y,t-\tau)\;dy d\tau.\]

 Consider, for $n=1,2,\dots$ and $\delta=T/n$, the mollified Cauchy problem
\begin{equation}
\left\{
\begin{aligned} \label{Approx} 
 &\partial_t \theta^n +\Psi_{\delta}(u^{n-1}) \cdot \nabla \theta^n =  \Delta \theta^n\\
 &\partial_t u^n +\nabla\cdot (\Psi_{\delta}(u^{n})\otimes u^n)+\nabla p^n= \Delta u^n+ \theta^{n} e_3\\
 &\nabla\cdot u^n=0.\\
\end{aligned}
\right.
\qquad x\in \R^3, t\in \R_{+}
\end{equation}
with data 
\bg
\label{mol-data}
\theta^n|_{t=0}=\theta_0
\qquad \hbox{and}\qquad u^n|_{t=0}=u_0. 
\ed
The iteration scheme starts with $u^0=0$.

Note that since $\hbox{div}\, u = 0$, we also have  $\mbox{div}\left(\Psi_{\delta}( u^n)\right) = 0$, for $t \in \R_{+}$.
At each step $n$, one solves recursively $n+1$ {\it linear\/} equations:
first one solves the transport-diffusion equation (with smooth convective velocity) for
the temperature; after $\theta^n$ is computed,
solving the second of~\eqref{Approx} amounts to solving a linear equation on each strip
$\R^3\times \hbox{$(m\delta,(m+1)\delta)$}$, for  $m=0,1,\ldots, n-1$.

For solutions to (\ref{Approx}) we have the following existence and uniqueness result,
\begin{proposition}
\label{pr:approx}
Let $(\theta_0,u_0) \in  L^2 \times L^2_{\boldsymbol{\sigma}} $. 
For each $n\in \{1,2,\dots\}$
 there exists a unique  weak solution $(\theta^n, u^n,p^n)$ of the approximating equations  with  data (\ref{mol-data}) such that, for any $T>0$,
$$\theta^n \in L^2(0,T; H^1) \cap L^{\infty}(0,T; L^2),\qquad
  u^n \in L^2(0,T; V) \cap L^{\infty}(0,T; L^2_{\boldsymbol{\sigma}})$$
and
$$ p^n\in L^{5/3}(0,T;L^{5/3})+L^\infty(0,T;L^6)$$
Moreover, 
for all $t>0$, $\theta^n$ and $u^n$ satisfy the energy inequalities as in~\eqref{eit} and \eqref{eiu}.
In particular, the sequences $\theta^n$, $u^n$ and $p^n$, $n=1,2,\dots$ are bounded in their respective spaces.
\end{proposition}

\pf

This can be proved using the Faedo-Galerkin method. As the the argument is standard
(see, {\it e.g.\/} \cite[Theorem~1.1, Chapter III]{Temam}, or \cite[Appendix]{CKN}),
we skip the details.
We only prove the condition on the pressure since this is the only change that we have to make in~\cite{CKN}.

 Taking the divergence of the second equation in (\ref{Approx}) we get $p=p_1^n+p_2^n$, where
\[\Delta p_1^n = -\sum_{i,j} \partial_i \partial_j (u_i^n u_j^n) \]
and
\[\Delta p_2^n = \partial_{x_3} \theta^n.\]
Thus, $p_1^n\in L^{5/3}(0,T;L^{5/3})$ uniformly with respect to~$n$, by the energy inequality for~$u^n$, interpolation,
and the Calderon-Zygmund theorem, as proved in~\cite{CKN}.
On the other hand, $-p_2^n=E_{x_3}*\theta$, where $E(x)$ is the fundamental solution of $-\Delta$.
Thus, $E_{x_3}(x)=\frac{c\,x_3}{|x|^3}$ belongs to the Lorentz space $L^{3/2,\infty}(\R^3)$.
But $\theta^n\in L^\infty(0,T;L^2)$
uniformly with respect to~$n$, hence
Young convolution inequality in Lorentz spaces
(see \cite[Chapter~2]{Lem02}) yields $p_2^n\in L^\infty(0,T;L^6)$
uniformly with respect to~$n$.

\cbdu

It follows from Proposition~\ref{pr:approx} that, extracting suitable subsequences,
$p^n=p_1^n+p_2^n$, where $(p_1^n)$ converges weakly in $L^{5/3}(0,T;L^{5/3})$ and
$p_2^n$ converges in $L^\infty(0,T;L^6)$ in the weak-* topology.
Moreover, $(u^n)$ is convergent with respect to the topologies 
listed in~\cite[p. 828-829]{CKN}.
On the other hand $(\theta^n)$ will be convergent with respect to the same topologies,
because all estimates available for $u^n$  hold also for~$\theta^n$.

No additional difficulty in the passage to the limit in the nonliner terms arises,
in the equation of the temperature  other those already existing
for the Navier--Stokes equations.
Hence, the distributional limit $(\theta,u,p)$ of a  convergent subsequence of $(\theta^n,u^n,p^n)$
is a weak solutions of the Boussinesq system.
This establishes Proposition~\ref{th:existen}.

\cbdu

\medskip
We finish this section by establishing the natural $L^p$-estimates for the approximating temperatures:

\begin{lemma}
 \label{lemma-lpd}
Let  $(\theta_0,u_0)\in L^2\times L^2_{\boldsymbol{\sigma}}$ and
let $(\theta^n,u^n,p^n)$ be the solution of the mollified Boussinesq system~\eqref{Approx} for some $n\in \{1,2,\ldots\}$.
Let also $1\le p<\infty$. If $\theta_0\in L^1\cap L^p$, then
\begin{equation}
\label{fo-de1} 
  \|\theta^n(t)\|_p \le \|\theta_0\|_1 
    \Bigl( \frac{c\,t}{p} + A\Bigr)^{-\frac{3}{2}(1-\frac{1}{p})},
\end{equation}
where $A=A(p,\|\theta_0\|_1,\|\theta_0\|_p)$ and $c>0$ is an absolute constant.
\end{lemma}

\pf
First notice that, for each $n$, $\theta^n$ is the solution of a linear transport-diffusion equation with smooth and 
divergence-free velocity $\Psi(u^{n-1})$.
The $L^p$ decay estimates for these equations are well known.
We reproduce the same proof as in~\cite{CC, EscZ91}) expliciting better the constants, as we will need the expressions
of such constants later on.

A basic estimate (valid for $1\le p\le \infty$) is
\begin{equation}
\label{L1b}
 \|\theta^n(t)\|_p \le \|\theta_0\|_p.
\end{equation}
See  \cite[Corollary 2.6]{CC} for a nice proof of~\eqref{L1b} that remains valid in the much more general case of
trasport equations with (or without) fractional diffusion.

We start with the case $2\le p<\infty$.
Multiplying the equation for $\theta^n$ by $p|\theta^n|^{p-2}\theta^n$ and integrating we get
\[
 \frac{d}{dt}\|\theta^n(t)\|_p^p +\frac{4(p-1)}{p} \|\nabla( |\theta^n|^{p/2} ) (t)\|^2 \le 0.
\]
By the Sobolev embedding theorem, $\dot H^1\subset L^6$, hence
\[
 \|\theta^n(t)\|_{3p}^p\le C\| \nabla ( |\theta^n|^{p/2} )(t) \|^2.
\]
The interpolation inequality yields
\[
 \|\theta^n(t)\|_p\le \|\theta^n(t)\|_1^{2/(3p-1)} \|\theta^n(t)\|_{3p}^{3(p-1)/(3p-1)}.
\]
Combining these two inequalities with the basic estimate $\|\theta^n(t)\|_1\le\|\theta_0\|_1$,
we obtain the differential inequality
\[
 \frac{d}{dt}\|\theta^n(t)\|_p^p \le -\frac{4(p-1)}{Cp\|\theta_0\|_1^{2p/(3p-3)}}  \bigl( \|\theta^n(t)\|_p^p \bigr)^{1+\frac{2}{3p-3}}.
\]
Integrating this we get
\[
 \|\theta^n(t)\|_p^p \le 
\biggl( \frac{ 8\,t}{3Cp \|\theta_0\|_1^{(2p)/(3p-3)} } \, +\, \frac{1}{\|\theta_0\|_p^{(2p)/(3p-3)} } \biggr)^{-3(p-1)/2}.
\]
and estimate~\eqref{fo-de1} follows with
$$ A=\left(\frac{\|\theta_0\|_1}{\|\theta_0\|_p}\right)^{(2p)/(3p-3)}.$$

The case $1\le p<2$ is deduced by interpolation.

\cbdu

\medskip
In the $p=2$ case, we obtain the following

\begin{lemma}
 \label{lemma-2da}
 Let $\theta_0\in L^1\cap L^2$ and $u_0\in L^2_{\boldsymbol{\sigma}}$.
Let $(\theta^n,u^n,p^n)$ be the solution of the mollified Boussinesq system~\eqref{Approx} for some $n\in\{1,2,\ldots\}$.
Then,
\begin{equation}
\label{theta-2} 
\begin{split}
&  \|\theta^n(t)\| \le \|\theta_0\|_1 
    \bigl( C\,t + A_0\bigr)^{-3/4},\\
& \|u^n(t)\| \leq \|u_0\|+ C'\|\theta_0\|_1\,t^{1/4},
\end{split}
\end{equation}
for two absolute constants $C,C'>0$ and $A_0=\bigl(\|\theta_0\|_1/\|\theta_0\|_2\bigr)^{4/3}$.
\end{lemma}

\pf
We only have to estimate the $L^2$ norm of the velocity.
We make use of the identity 
$$\frac{d}{dt} \|u^n(t)\|^2_2 = \int u^n \partial_t u^n,$$
that can be justified exactly as for the mollified Navier-Stokes equations, see \cite{CKN} and \cite{Schon84}.

Multiplying the velocity equation in the Boussinesq system (\ref{Approx})  by $u^n$ and integrating
we get
\begin{equation}
 \label{cosa1}
\frac{d}{dt}\|u^n(t)\|^2 + 2\|\nabla u^n(t)\|^2 \leq 2  \|u^n(t)\|\|\theta^n(t)\|. 
\end{equation}
Dividing by $ \|u^n(t)\|$,
\[
 \frac{d}{dt}\|u^n(t)\| \leq   \|\theta^n(t)\|.
\]
Now we use the decay of $\|\theta^n(t)\|$ obtained in Lemma~\ref{lemma-lpd}.
Integrating we obtain the second of~\eqref{theta-2}.

\cbdu

\section{Improved bounds for weak solution in the case $\int\theta_0=0$}
\label{secsix}

 The estimates obtained in Lemma~\ref{lemma-2da} can be considerably improved provided we
 additionally assume $\int\theta_0=0$ and the moment condition $\theta_0\in L^1_1$.
First of all, from an elementary heat kernel estimate one easily checks that in this case
\begin{equation}
 \label{hke}
\|e^{t\Delta}\theta_0\|^2 \le A_1(t+1)^{-5/2},
\end{equation}
where $A_1>0$ depends only on the data, through its $L^2$-norm and  the $\int |x|\,|\theta_0(x)|\,dx$ integral.

We will now see that in this case the approximated 
temperature $\theta^n(t)$ also decays at the faster rate $(t+1)^{-5/4}$ in the $L^2$-norm.
 Using this new decay rate for $\theta^n$ it is possible to show that the velocities
 are uniformly  bounded in~$L^2$.
Once we have a solution such that $\|u^n(t)\|$ remains bounded as $t\to\infty$ 
one can go further and prove that $u^n$ actually decays at infinity in $L^2$ at some
algebraic decay rate, depending on the decay of the linear evolution~$e^{t\Delta}u_0$.
In view of the passage to the limit from the mollified system~\eqref{Approx} to the Boussinesq system~\eqref{B},
all the estimates must be independent on~$n$.

 \begin{proposition}\label{th:firstdecay}
Let $(\theta_0,u_0) \in   (L^1_1\cap L^2)\times L^2_{\boldsymbol{\sigma}} $ and assume $\int\theta_0=0$.
There exists an absolute constant $\varepsilon_0>0$ such that if
\begin{equation}
 \label{vares}
\|\theta_0\|_1<\varepsilon_0
\end{equation}
then the solution of  the mollified Boussinesq system (\ref{Approx}), with data $(u_0,\theta_0) $ satisfies
 \begin{align}
\|\theta^n(t)\| ^2\le& A(1+t)^{-\frac{5}{2}}
\end{align}
and
\begin{align}
\|u^{n}(t)\|^2 \le A
\end{align}
for all $n \in \N$ and $t\in\R^+$.
Here $A>0$ is some constant depending  on the data $u_0$ and  $\theta_0$, and independent on~$n$, and $t$.
 \end{proposition}

\pf  
We denote by $C$ a positive absolute constant, which may change from line to line.
We also denote by $A_1, A_2,\ldots$ positive constants that depend only on the data.
More precisely, $A_j=A_j \bigl(\|\theta_0\|,\|\theta_0\|_{L^1_1},\|u_0\|\bigr)$.

\medskip
{\bf Step 1:} An auxiliary estimate.

\medskip

We make use of the Fourier splitting method introduced in~\cite{Schon84}.
The first step consists in multiplying the temperature equation  by $\theta^n$ and to integrate by parts.
Using the Plancherel theorem in the energy inequality for $\theta^n$, we get
 \begin{equation*} 
\label{FS2}
\frac{1}{2}\frac{d}{dt} \int |\widehat{\theta^n}(\xi,t)|^2 \; d\xi \leq -\int |\xi|^2|\widehat{\theta^n}(\xi,t)|^2 \; d \xi
\end{equation*}
 Now split the integral on the right hand side  into $\Ss \cup \Ss^c$,
where
\begin{equation} \label{def-S}
\Ss = \left\{ \xi: |\xi| \le\left(\frac{k}{2(t+1)}\right)^{1/2}\right\},
\end{equation}
and $k$ is a constant to be determined below.
Noting that for $\xi  \in \Ss^c$
one has
$-|\xi |^{2} \leq - \frac{k}{2(t+1)}$, it follows that
\begin{equation*}
\begin{split}
\frac{d}{dt}\int |\widehat{\theta^n}(\xi,t)|^2 \;d\xi
&\leq - \frac{k}{t+1}\int_{\Ss^c}|\widehat{\theta^n}(\xi,t)|^2 \\
&= - \frac{k}{t+1}\int |\widehat{\theta^n}(\xi,t)|^2\;d\xi+
\frac{k}{t+1}\int_{\Ss}|\widehat{\theta^n}(\xi,t)|^2\;d\xi.
\end{split}
\end{equation*}
Multiplying by $(1+t)^k$ we obtain
\begin{equation} \label{FS4}
\frac{d}{dt}\left[(1+t)^{k}\int |\widehat{\theta^n}(\xi,t)|^2 \;d\xi\right]
\leq k(t+1)^{k-1}\int_{\Ss} |\widehat{\theta^n}(\xi,t)|^2 d\xi,
\end{equation}
where $\Ss$ is as in \eqref{def-S}.

Taking the Fourier transform in the equation for $\theta^n$ in~\eqref{Approx} we get
 \begin{equation}
 \label{FS-theta}|\widehat{\theta^n}(\xi,t)|^2
 \leq 2|e^{-t|\xi|^2}\widehat{\theta_0}|^2 +2|\xi|^2\left(\int_0^t  \|u^{n-1}(s)\| \, \|\theta^n(s)\| ds\right)^2.
\end{equation}

Hence 
\[
\int_{\Ss} |\widehat{\theta^n}(\xi,t)|^2 \;d\xi \leq C \left[\|e^{-|\xi|^2t}\theta_0\|^2 
  + (1+t)^{-5/2}\left(  \int_0^t\|u^{n-1}(s)\| \,\|\theta^n(s)\|\, ds  \right)^2 \right].
\]

Replacing this in~\eqref{FS4} and applying the Plancherel theorem we get
\[\frac{d}{dt}\left[(1+t)^{k} \|\theta^n(t)\|^2\right] 
 \leq C\left[\|e^{t\Delta}\theta_0\|^2(1+t)^{k-1} + (1+t)^{k-7/2}\left(  \int_0^t\|u^{n-1}\|\, \|\theta^n\| ds  \right)^2\right]. \]

From where it follows, letting $k=7/2$,
\begin{equation}  
\label{ bueno}
\begin{split}
\|\theta^n(t)\|^2 \leq (1+t)^{-7/2}
\biggl\{ \|\theta_0\|^2  &+ C\biggl[ \biggl( \int_0^t \|e^{s\Delta}\theta_0\|^2(1+s)^{5/2}ds\biggr)\\
  &\qquad \quad+\int_0^t \left(\int_0^s \|u^{n-1}(r)\| \, \|\theta^n(r)\|dr\right)^2\, ds   \biggr]  \biggr\}.
\end{split}
 \end{equation}

Recalling the estimate~\eqref{hke} for the linear evolution we get, for all $n \in \N$,
\begin{equation}
\label{bueno1}
\|\theta^n(t)\|^2 
\leq C\left[A_1(1+t)^{-5/2}  + (1+t)^{-7/2}\int_0^t \left(\int_0^s \|u^{n-1}(r)\| \,\|\theta^n(r)\|dr\right)^2 ds  .  \right]
\end{equation}

We now use the following inequality, deduced from estimate~\eqref{theta-2},
\begin{equation}
\label{moreco} 
\|\theta^n(t)\|\le C\,\|\theta_0\|_1\,t^{-3/4}.
\end{equation}

Putting this inside~\eqref{bueno1} we obtain a new bound for~$\|\theta^n\|^2$, namely
\begin{equation}
\label{buenoo1}
\|\theta^n(t)\|^2 \leq C\left[ A_1(1+t)^{-5/2}  + \|\theta_0\|_1^2(1+t)^{-7/2}\int_0^t \left(\int_0^s \|u^{n-1}(r)\| \, r^{-3/4}\,dr\right)^2 ds .   \right]
\end{equation}

\medskip
{\bf Step 2:}  The inductive argument.

\medskip
We will now prove by  induction that, for all positive integer $n$ we have
\begin{align}\label{HI}
\|u^{n-1}(t)\| \leq \|u_0\|+M t^{1/8}
\end{align}
where $M>0$ is some constant independent on~$n$ (but possibly dependent on the data $\theta_0$,  $u_0$)
to be determined.
Notice that estimate~\eqref{HI} is actually better than what we have so far (compare with the second of~\eqref{theta-2}).

For $n=1$ the inductive condition~\eqref{HI} is immediate since $u^0=0$. 
Let us now prove that
$\|u^{n}\| \leq \|u_0\|+M t^{1/8}$ assuming that~\eqref{HI} holds true.

\medskip
We get from~\eqref{buenoo1} and the induction assumption~\eqref{HI}
\begin{equation*}
\|\theta^n(t)\|^2\le C\biggl[ A_1(1+t)^{-5/2}+\|\theta_0\|_1^2\|u_0\|^2(1+t)^{-2}+M^2\|\theta_0\|_1^2(1+t)^{-7/4}\biggr].
\end{equation*}
This implies
\begin{equation}
\label{buenooo}
\|\theta^n(t)\|\le C\biggl[ \sqrt A_1+\|\theta_0\|_1\Bigl(\|u_0\|+M\Bigr)\biggr](1+t)^{-7/8}.
\end{equation}

This last inequality  will be used to estimate $\|u^n\|$ as follows. 
First recall that
 \begin{equation}
\label{difino}
\frac{d}{dt} \|u^n(t)\| \leq \|\theta^n(t)\|.
\end{equation}
After an integration in time we get, using \eqref{buenooo},
\[
 \|u^n(t)\|\le \|u_0\|+C\biggl[\sqrt A_1+\|\theta_0\|_1\Bigl(\|u_0\|+M\Bigr)\biggr]t^{1/8}.
\]
For the induction argument we need to prove $\|u^n\|\le \|u_0\|+Mt^{1/8}$. 
Hence we need that
\[
 C\biggl[\sqrt A_1+\|\theta_0\|_1\Bigl(\|u_0\|+M\Bigr)\biggr]\le M.
\]
Choosing $M$ large enough, for example,
$$M=\max\{\|u_0\|, 2C\sqrt A_1\},$$
our condition then boils down to the inequality
\[
 \|\theta_0\|_1\le1/(4C).
\]
The validity this last inequality is insured by assumption~\eqref{vares}.
This concludes the induction argument and establishes the validity of the estimate~\eqref{HI} for all~$n$.

\medskip
{\bf Step 3:} Uniform bound for the $L^2$-norm of the velocities $u^n$.

\medskip
The result of Step 2 implies the existence of a constant $A_2>0$ such that, for all $n\ge1$,
\begin{equation} 
 \label{tbco}
\|u^{n-1}(t)\|^2\le A_2(1+t)^{1/4}.
\end{equation}
In the proof of the previous Step we also deduced that, for some $A_3>0$,
$$ \|\theta^n(t)\|^2\le A_3(1+t)^{-7/4}.$$
Combining such two estimates  with inequality~\eqref{bueno1} we easily get
\[
\label{e-2}
 \|\theta^n(t)\|^2\le A_4(1+t)^{-2}.
\]
Now using  this improved estimate for~$\|\theta^n\|^2$ with~\eqref{tbco} in~\eqref{bueno1}
arrive at
\[
\label{e-94}
 \|\theta^n(t)\|^2\le A_5(1+t)^{-9/4}.
\]
Going back to the differential inequality~\eqref{difino} we finally get, for some constant $A_6>0$ independent on~$n$, and $t\in\R^+$,
\[
 \|u^n(t)\|^2\le A_6.
\]
Replacing in~\eqref{bueno1} we can further improve the decay of~$\theta^n$ up to
\[
\label{e-52}
 \|\theta^n(t)\|^2\le A_7(1+t)^{-5/2}.
\]

\endProof

\medskip

We now would like to improve the result of the previous Proposition by establishing decay properties for $\|u^n(t)\|^2$.
Specifically, if we assume in addition that the linear part of the velocity satisfies 
$\|e^{t\Delta}u_0\|^2\le C(1+t)^{-1/2}$ (this happens {\it e.g.\/} when $u_0\in L^{3/2}\cap L^2_{\boldsymbol{\sigma}}$),
then the same decay holds for the approximate velocities  $u^n$.

\begin{proposition} \label{th:u-decay1}
 Let $(\theta_0,u_0) \in   (L^1_1\cap L^2 )\times L^2_{\boldsymbol{\sigma}})$.
 Assume also that $\int\theta_0=0$ and   $\|\theta_0\|_1<\varepsilon_0$,
where $\varepsilon_0$ is the constant obtained in the previous proposition.

Then the approximate solutions of~\eqref{Approx} satisfy,
\begin{equation*}
 \| u^n(t)\|\to 0, \qquad\hbox{as $t\to\infty$},
\end{equation*}
uniformly with respect to~$n$.
Moreover, if $u_0\in L^{3/2}\cap L^2_{\boldsymbol{\sigma}}$, then
 \begin{equation}
\|u^n(t)\| ^2\le A(1+t)^{-1/2},
\end{equation}
for some constant $A>0$ independent on $n$, and $t$.
\end{proposition} 

\pf 
We denote by $A>0$ a constant depending only on the data that might change from line to line.
The proof follows by Fourier Splitting. 
 Since the estimates are independent of $n$,  we simply denote the solutions by $(\theta,u)$.
Multiply the second equation in (\ref{Approx}) by $u$, integrate in space. 
By~Proposition~\ref{th:firstdecay} we get
\begin{equation}
 \label{u-decay0} \frac{d}{dt}  \|u(t)\|^2  + 2  \|\nabla u\|^2  \leq  C(t+1)^{-\frac{5}{4}}
\end{equation}

Arguing as for the proof of inequality~\eqref{FS4} we obtain
\begin{equation}
 \label{u-decay}
\frac{d}{dt}\biggl[(t+1)^k \|u(t)\|^2\biggr] \leq Ck(t+1)^{k-1} \int_{\Ss} |\widehat{u}(\xi,t)|^2 \;d\xi +C (t+1)^{k-\frac{5}{4}}
\end{equation}
Where $\Ss$ was defined in (\ref{def-S}). From now on, $k=7/2$.

We need to estimate $|\widehat{u}(\xi,t)| $ for $\xi \in \Ss$. Computing the Fourier transform
in the equation for~$u$ in~\eqref{Approx}, next applying  the estimate $\|u(t)\|^2\le A$ obtained in
Proposition~\ref{th:firstdecay} we get
\[
\begin{split}
\label{ufourier} |\widehat{u}(\xi,t)| 
&\leq e^{-t|\xi|^2} |\widehat{u_0}| +|\xi|\int_0^t\|u(s)\|^2\,ds + \int_0^t |\widehat\theta(\xi,s)|\,ds\\
&\leq e^{-t|\xi|^2} |\widehat{u_0}| +At|\xi| + \int_0^t |\widehat\theta(\xi,s)|\,ds.
\end{split}
\]
But computing the Fourier transform
in the equation for~$\theta$ in~\eqref{Approx} and applying  once more the estimates of Proposition~\ref{th:firstdecay}
 we have
\[
|\widehat \theta(\xi,s)|\le  e^{-s|\xi|^2} |\widehat\theta_0| +A|\xi|
\le|\widehat\theta_0(\xi)| +A|\xi|
\]
Hence,
\begin{equation}
\label{taks}
|\widehat{u}(\xi,t)|^2 \leq 
A\Bigl[ e^{-2t|\xi|^2} \bigl( |\widehat{u_0}|^2 + |\widehat{\theta_0}|^2 \bigr)     +t^2|\xi|^2 \Bigr].
\end{equation} 
Integrating on~$\Ss$ and applying inequality~\eqref{hke} we deduce
\[
 \int_{\Ss} |\widehat u(\xi,t)|^2\,d\xi\le A\biggl[ \|e^{t\Delta}u_0\|^2+ (1+t)^{-1/2}\biggr].
\]
Putting this inside~\eqref{u-decay} and integrating on an interval of the form
$[t_\epsilon,t]$,
with $\epsilon>0$ arbitrary and $t_\epsilon$ chosen
in a such way $\|e^{t\Delta}u_0\|^2<\epsilon$ for $t\ge t_\epsilon$,
we obtain
$$\|u(t)\|\to0,\qquad \hbox{as $t\to\infty$}.$$

On the other hand, in the case $u_0\in L^{3/2}\cap L^2_{\boldsymbol{\sigma}}$,
we have
\begin{equation*} 
 \int_{\Ss} |\widehat u(\xi,t)|^2\,d\xi\le A(1+t)^{-1/2}.
\end{equation*}
Now going back to~\eqref{u-decay} and integrating in time we finally get
\begin{equation*}
 \|u(t)\|^2\le A(1+t)^{-1/2}.
\end{equation*}

\endProof

\bigskip

We are now in the position of deducing our result on weak solutions to the Boussinesq system~(\ref{B}).

\medskip
\paragraph{\bf Proof of Theorem~\ref{th:belo}.}

Now this is immediate:
passing to a subsequence, the approximate solutions $\theta^n$ and $u^n$ 
converge in $L^2_{loc}(\R^+,\R^3)$ to a weak solution $(\theta,u)$ of the Boussinesq system~\eqref{B}.
Moreover, the previous Lemmata imply that $\theta^n$ and $u^n$ satisfy  estimates of the form
$$ \|v^n(t)\| \le f(t), \qquad \hbox{for all $t>0$},$$
where $f(t)$ is a continuous function independent on~$n$.
Then the same estimate must hold for the  limit $\theta$ and $u$, except possibly points in a set of measure zero.
But since weak solutions are necessarily continuous from $[0,\infty)$ to $L^2$  under the weak topology,
$\|\theta(t)\|$ and $\|u(t)\|$ are lower semi-continuous and hence they satisfy the above estimate for all~$t>0$.
This observation on the weak semi-continuity is borrowed from~\cite{KM}.
\cbdu

\medskip
\section{Strong solutions: preliminary lemmata}
\label{sec:seven}

The  integral formulation  for the Boussinesq system, formally equivalent to~\eqref{B} 
reads
\begin{equation}
 \label{IE}
\left\{
\begin{aligned}
 &\theta(t)=e^{t\Delta}\theta_0-\int_0^t e^{(t-s)\Delta}\nabla \cdot(\theta u)(s)\,ds\\
 &u(t)=e^{t\Delta}u_0-\int_0^t e^{(t-s)\Delta}\P\nabla \cdot(u\otimes u)(s)\,ds
+\int_0^t e^{(t-s)\Delta}\P \theta(s) e_3\,ds.\\
&\nabla \cdot u_0=0
\end{aligned}
\right.
\end{equation}

The above system will be solved applying the following abstract lemma, which slightly generalizes that
of  G.~Karch et N.~Prioux  (see  \cite[Lemma~2.1]{KarPri08}).

\begin{lemma}
\label{LKP}
Let $\mathcal{X}$ and $\mathcal{Y}$ be two Banach spaces, let
$B\colon \mathcal{X}\times\mathcal{X}\to\mathcal{X}$
and
$\widetilde B\colon\mathcal{Y}\times\mathcal{X}\to\mathcal{Y}$
be two bilinear maps
and $L\colon\mathcal{Y}\to\mathcal{X}$ a linear map
satisfying the estimates
$\|B(u,v)\|_{\mathcal{X}}\le \alpha_1\|u\|_{\mathcal{X}}\|v\|_{\mathcal{X}}$,
$\|\widetilde B(\theta,v)\|_{\mathcal{Y}}\le \alpha_2\|\theta\|_{\mathcal{Y}}\|u\|_{\mathcal{X}}$
and
$\|L(\theta)\|_{\mathcal{X}}\le \alpha_3\|\theta\|_{\mathcal{Y}}$,
for some positive constants $\alpha_1$, $\alpha_2$ and  $\alpha_3$.

Let $0<\eta<1$ be arbitrary.
For every $(U,\Theta)\in \mathcal{X}\times\mathcal{Y}$ such that
\begin{equation*}
\eta\|U\|_{\mathcal{X}}+\alpha_3\|\Theta\|_{\mathcal{Y}}
\le \frac{\alpha_1\eta(1-\eta)^2}{(2\alpha_1+\alpha_2)^2}, 
\end{equation*}
the system
\begin{equation}
\label{ABS}
 \theta=\Theta+\widetilde B(\theta,u), \qquad u=U+B(u,u)+L(\theta)
\end{equation}
has a solution $(u,\theta)\in\mathcal{X}\times\mathcal{Y}$.
This is the unique solution satisfying the condition
\begin{equation*}
\eta\|u\|_\mathcal{X}+\alpha_3\|\theta\|_\mathcal{Y}\le \eta(1-\eta)(2\alpha_1+\alpha_2).
\end{equation*}
\end{lemma}
\Proof
In the case $0<\alpha_3<1$, one can take $\eta=\alpha_3$.
In such particular case, this lemma is already known, see~\cite[Lemma~2.1]{KarPri08}.
Therefore, we only have to prove that we can get rid of the restriction $0<\alpha_3<1$.
This is straightforward. We introduce on the space $\Yy$ an equivalent norm,
defined by $\|\theta\|_{\Yy'}=\frac{\alpha_3}{\eta}\|\theta\|_{\Yy}$.
By Lemma~2.1 of Karch and Prioux, applied in the space $(\Xx,\|\cdot\|_\Xx)$ and $(\Yy,\|\cdot\|_{\Yy'})$,
we have that if
\begin{equation*}
\|U\|_{\Xx}+\|\Theta\|_{\Yy'}
\le \frac{\alpha_1(1-\eta)^2}{(2\alpha_1+\alpha_2)^2}, 
\end{equation*}
then the system~\eqref{ABS} has a unique solution such that
\begin{equation*}
\|u\|_\mathcal{X}+\|\theta\|_\mathcal{Y'}\le (1-\eta)(2\alpha_1+\alpha_2).
\end{equation*}
The conclusion of Lemma~\ref{LKP} is now immediate.

\endProof

\begin{remark}
Using our improved version of Lemma~2.1 of \cite{KarPri08}, it is be  possible to get rid of 
the smallness assumption $|\beta|<1$ in the main results of Karch and Prioux~\cite{KarPri08}.
\end{remark}

\begin{remark}
The proof of Lemma~2.1 in~\cite{KarPri08} relies on the 
contraction mapping theorem. 
In particular, the solution can be obtained passing to the limit with respect to the $\mathcal{X}\times\mathcal{Y}$-norm
in the iteration scheme $(k=1,2,\ldots)$
\begin{equation}
\label{ites}
\begin{split}
&(u^0,\theta^0)=(U,\Theta),\\
&(u^{k+1},\theta^{k+1})=\bigl(u^0+B(u^k,u^k)+L(\theta^k),\theta^0+\widetilde B(u^k,\theta^k)\bigr).
\end{split}
\end{equation}
See~\cite{KarPri08} for more details. See also~\cite[Lemma~4.3]{Prioux07} for similar abstract lemmata.
\end{remark}

Let $a\ge1$. We define $\Xx_a$ as the Banach space of divergence-free vector fields $u=u(x,t)$,
defined and measurable on $\R^3\times\R^+$,
such that, for some $C>0$,
\begin{equation}
\label{edeca1}
|u(x,t)|\le C\inf_{0\le \eta\le a} |x|^{-\eta}(1+t)^{(\eta-1)/2}.
\end{equation}
In the same way, for $b\ge 3$ we define the space $\Yy_b$ of functons $\theta\in L^\infty_t(L^1)$ satisfying the estimates
\begin{equation}
\label{edeca11}
|\theta(x,t)|\le C\inf_{0\le \eta\le b} |x|^{-\eta}(1+t)^{(\eta-3)/2}.
\end{equation}
Such spaces are equipped with their natural norms.

 They are obviously decreasing with respect to inclusion as $a$ and $b$ grow.
Recalling the definition of~$\Xx$ and $\Yy$ in Section~\ref{sec:two.one},
we see that $\Xx_1=\Xx\cap L^\infty_{x,t}$ with equivalence of the norms
and that $\Yy_3=\Yy\cap L^\infty_{x,t}$.

The estimates~\eqref{deca1}-\eqref{deca11} in Proposition~\ref{prop2}
are thus equivalent to the conditions~$u\in \Xx_a$ and $\theta\in\Yy_b$.

\bigskip
We start with some elementary embeddings.

\begin{lemma}
 \label{emb}
Let $L^{p,q}$ be the Lorentz space,  with $1<p<\infty$ and $1\le q\le\infty$.
Then the following four inequalities hold:
\begin{equation}
 \label{elx}
\begin{split}
&\|u(t)\|_{L^{p,q}}\le C\|u\|_{\Xx}\,t^{\frac{1}{2}(\frac{3}{p}-1)}, \qquad 3<p\le \infty,\\
&\|u(t)\|_{L^{p,q}}\le C\|u\|_{\Xx_a}
 (1+t)^{\frac{1}{2}(\frac{3}{p}-1)}, \qquad \textstyle\frac{3}{a}<p\le\infty,\\
&\|\theta(t)\|_{L^{p,q}}\le C\|\theta\|_{\Yy}\, t^{\frac{1}{2}(\frac{3}{p}-3)},
\qquad 1<p\le\infty\\
&\|\theta(t)\|_{L^{p,q}}\le C\|\theta\|_{\Yy_3} (1+t)^{\frac{1}{2}(\frac{3}{p}-3)},
\qquad 1<p\le\infty\\
\end{split}
\end{equation}
for some constants $C$ depending only on $p$ and $q$.
In particular, choosing $p=q$ one gets the correponding estimates for the classical $L^p$-spaces.
\end{lemma}

\Proof
The above estimates for the weak-Lebesgue spaces $L^{p,\infty}$ are simple.
Indeed,  if $u\in\Xx$, then $|u(x,t)|\le C|x|^{-3/p}t^{\frac{1}{2}(\frac{3}{p}-1)}$ and 
one  has only to recall that any function bounded by  $|x|^{-3/p}$ belongs to $L^{p,\infty}$.
The other $L^{p,\infty}$-estimates  for $u$ and $\theta$ 
contained in~\eqref{elx} 
In the case $1\le q<\infty$, we use that $L^{p,q}$  is a real interpolation space
between $L^{p-\varepsilon,\infty}$ and $L^{p+\varepsilon,\infty}$, 
for all $1<p-\varepsilon<p<p+\varepsilon<\infty$.
Therefore estimates~\eqref{elx} for all $1\le q\le\infty$ follow from the corresponding estimates
in the particular case $q=\infty$ {\it via\/} the interpolation inequality.

\endProof

\bigskip
The first useful estimate in view of the application of Lemma~\ref{LKP} is the following.

\begin{lemma}
\label{L34}
Let $1\le a<3$ and $\theta\in \Yy_3$. We have, for some~constant $C>0$ depending only on $a$,
\begin{equation}
 \label{ltx}
\|L(\theta)\|_{\Xx_a}\le C\|\theta\|_{\Yy_3}.
\end{equation}
Moreover,
\begin{equation}
\label{obm}
\| L(\theta)\|_{\Xx}\le C\|\theta\|_{\Yy}.
\end{equation}
\end{lemma}

\Proof
We prove only~\eqref{ltx} since 
the proof of~\eqref{obm} is essentially the same.
By a renormalization, we can and do assume that $\|\theta\|_{\Yy_3}=1$.
Let $\K(x,t)$ be the kernel of the operator~$e^{t\Delta}\P$.
Then we can write
\begin{equation*}
L(\theta)(x,t)=\int_0^t\!\!\int \K(x-y,t-s)\theta(y,s)e_3\,dy\,ds.
\end{equation*}
We have the well known estimates for $\K$ (see, {\it e.g.\/}, \cite[Prop. 1]{Bra08})
\begin{equation}
 \label{Kest}
|\K(x,t)|\le C|x|^{-a}t^{-(3-a)/2}, \qquad\hbox{for  all $0\le a\le 3$}
\end{equation}
where $C>0$ is come constant independent on $x,t$ and $0\le a\le 3$.
We also recall the scaling relation
\begin{equation}
\label{scal}
 \K(x,t)=t^{-3/2}\K(x/\sqrt t,1)
\end{equation}
and the fact that $\K(\cdot,t)\in C^\infty(\R^3)$ for $t>0$.
The usual $L^p$ estimates for $\K$ are
\begin{equation}
 \label{lpek}
\|\K(t)\|_p\le Ct^{-3/2+\frac{3}{2p}}, \qquad 1<p\le \infty.
\end{equation}

Using the $L^2$-$L^2$ convolution inequality, we get
\begin{equation*}
\|L(\theta)(t)\|_\infty\le C\int_0^t(t-s)^{-3/4}(1+s)^{-3/4}\,ds\le C(1+t)^{-1/2}.
\end{equation*}
Owing to this estimate, the conclusion $L(\theta)\in \Xx_a$ will follow provided we prove the pointwise inequality,
\begin{equation*}
 |L(\theta)|(x,t)\le C_a|x|^{-a}t^{(a-1)/2}, \qquad\forall\,(x,t)\;\hbox{ s.t. }\; |x|\ge 2 \sqrt t.
\end{equation*}
This leads us to decompose
\begin{equation*}
 L(\theta)=I_1+I_2+I_3,
\end{equation*}
where $I_1=\int_0^t\!\!\int_{|y|\le |x|/2}\dots$, $I_2=\int_0^t\!\!\int_{|x-y|\le |x|/2}\dots$
and $I_3=\int_0^t\!\!\int_{|y|\ge |x|/2,\; |x-y|\ge |x|/2}\dots$.
Using $\theta\in L^\infty_t(L^1)$
we get
\begin{equation}
|I_1|(x,t)\le C|x|^{-3}t,
\end{equation}
which is even better in the region $\{(x,t)\colon |x|\ge 2\sqrt t\}$ than what we need (recall that $1\le a<3$).
Using now $|\theta(x,t)|\le C|x|^{-3}$  and the scaling properties of~$\K$ we obtain
by a change of variables
\begin{equation*}
 \begin{split}
  |I_2|(x,t) &\le C|x|^{-3} \int_0^t\!\!\int_{|y|\le |x|/(2\sqrt s)} |\K(y,1)|\,dy\\
& \le  C|x|^{-3}\,t\,\log(|x|/\sqrt t)\\
&\le \textstyle\frac{C}{3-a} |x|^{-a}t^{(a-1)/2}
\end{split}
\end{equation*}
for $|x|\ge 2\sqrt t$, and $1\le a<3$.
Next, using again $|\theta(x,t)|\le C|x|^{-3}$ and $|\K(x,t)|\le C|x|^{-3}$,
\begin{equation*}
 |I_3|(x,t) \le C\int_0^t\!\!\int_{|y|\ge |x|/2}|y|^{-6}\,dy\,ds\le C|x|^{-3}\,t.
\end{equation*}
Therefore,
$$ |L(\theta)|(x,t) \le \textstyle\frac{C}{3-a}|x|^{-a}t^{(a-1)/2},
\qquad |x|\ge 2\sqrt t.$$

%
%
Lemma~\ref{L34} in now established.

\endProof

\bigskip
We collect in the following Lemma all the estimates  on $B(u,v)$
that we shall need. (We will apply estimate~\eqref{bbe3} in the proof of Proposition~\ref{prop1},
estimate~\eqref{bbe2} for Proposition~\eqref{prop1} and estimate~\eqref{bbe2} in Theorem~\ref{theor1}.

\begin{lemma}
 \label{buv}
Let $1\le a<3$.
For some constant $C>0$, depending only on $a$ we have
\begin{equation}
\label{bbe1} 
\bigl\| B(u,v)\bigr\|_{\Xx_a}\le C\bigl\|u\bigr\|_{\Xx}\bigl\|v\bigr\|_{\Xx_a}\,\,
\end{equation}
and
\begin{equation}
\label{bbe2} 
\bigl\| B(u,v)\bigr\|_{\Xx_{(2a)_*}}\le C\bigl\|u\bigr\|_{\Xx_a}\bigl\|v\bigr\|_{\Xx_a},\,\,
\end{equation}
where $(2a)_*=\min\{2a,4\}$.
Moreover, 
\begin{equation}
\label{bbe3} 
\bigl\| B(u,v)\bigr\|_{\Xx}\le C\bigl\|u\bigr\|_{\Xx}\bigl\|v\bigr\|_{\Xx}\,\,.
\end{equation}
\end{lemma}

\Proof
We begin with the proof of  the first estimate.
As before, we can assume $\bigl\|u\bigr\|_{\Xx}=\bigl\|v\bigr\|_{\Xx_a}=1$.
We start writing
\begin{equation}
 \label{bwf}
B(u,v)(x,t)=\int_0^t\!\!\int F(x-y,t-s)(u\otimes v)(y,s)\,dy\,ds,
\end{equation}
where $F(x,t)$ is the kernel of the operator $e^{t\Delta}\P\nabla$.

The well known counterpart of relations~\eqref{Kest}-\eqref{scal} are
(see, {\it e.g.\/}, \cite{Miy00}, \cite[Prop. 1]{Bra08})
\begin{equation}
 \label{Fest}
|F(x,t)|\le C|x|^{-\eta}t^{-(4-\eta)/2}, \qquad\hbox{for all $0\le \eta\le 4$}
\end{equation}
and some constant $C>0$ independent on~$x,t$ and on $0\le a\le 4$.
Moreover,
\begin{equation}
\label{scalf}
 F(x,t)=t^{-2}F(x/\sqrt t,1).
\end{equation}
These bounds imply the useful
estimates
\begin{equation}
\label{l1f} 
\|F(t)\|_p\le Ct^{-2+\frac{3}{2p}} \qquad (1\le p\le \infty).
\end{equation}
Applying the first of~\eqref{elx} with $p=q=6$, we get $\|u(t)\|_6\le t^{-1/4}$.
Similarily,  $\|v(s)\|_\infty\le (1+t)^{-1/2}$. Hence,
\begin{equation}
\begin{split}
\label{libu}
\|B(u,v)\|_\infty
&\le C\int_0^t\|F(t-s)\|_{6/5}\|u\otimes v(s)\|_{6}\,ds\\
&\le C\int_0^t(t-s)^{-3/4}s^{-1/4}(1+s)^{-1/2}\,ds\\
&\le C(1+t)^{-1/2}.
\end{split}
 \end{equation}

It remains to establish a pointwise estimate in the in the region $\{(x,t)\colon |x|\ge 2\sqrt t$\}.
Let us decompose
$$ B(u,v)=I'_1+I'_2,$$
by splitting the integrals as $\int_0^t\!\!\int_{|y|\le |x|/2}\dots$ and $\int_0^t\int_{|y|\ge |x|/2}\dots$\,.
For the estimate of $I'_1$ we use $|u|\le s^{-1/2}$, $|v|\le |y|^{-a}s^{(a-1)/2}$ and $|F(x,t)|\le C|x|^{-3}t^{-1/2}$.
For the estimate of $I'_2$ we use again $|u|\le s^{-1/2}$, $|v|\le C|y|^{-a}s^{(a-1)/2}$, and the $L^1$-estimate
for $F$ (see ~\eqref{l1f}).
This leads to
\begin{equation*}
|B(u,v)|(x,t)\le Ct|x|^{-a}t^{(a-1)/2}.
\end{equation*}
The conclusion follows combining this with estimate~\eqref{libu}.
The proof of estimate~\eqref{bbe2} is similar. Notice the limitation $(2a)_*\le 4$, which is due to
the restriction on~$\eta$ in inequality~\eqref{Fest}.
The proof of~\eqref{bbe3} also follows along the same lines and is left to the reader.

\endProof

\bigskip
We finish with $\widetilde B(\theta,u)$.

\begin{lemma}
\label{buvv}
Let $a\ge 1$, $b\ge 3$. For some constant $C>0$ depending only on $a,b$, we have
\begin{equation}
\label{bte1} 
\bigl\| \widetilde B(\theta,u)\bigr\|_{\Yy_b}
\le C\bigl\|u\bigr\|_{\Xx}\bigl\|\theta\bigr\|_{\Yy_b}\,\,
\end{equation}
and
\begin{equation}
\label{bte2} 
\bigl\| \widetilde B(\theta,u)\bigr\|_{\Yy_{a+b}}
\le C\bigl\|u\bigr\|_{\Xx_a}\bigl\|\theta\bigr\|_{\Yy_b}\,\,.
\end{equation}

Moreover,
\begin{equation}
\label{bte3} 
\bigl\| \widetilde B(\theta,u)\bigr\|_{\Yy}
\le C\bigl\|u\bigr\|_{\Xx}\bigl\|\theta\bigr\|_{\Yy}\,\,.
\end{equation}
\end{lemma}

\Proof
As before, we give details only for the first estimate.
Denoting $\widetilde F(x,t)$ the kernel of $e^{t\Delta}\nabla$,
we can write
\begin{equation}
 \label{btwf}
\widetilde B(\theta,u)(x,t)=\int_0^t\!\!\int \widetilde F(x-y,t-s)(\theta\, u)(y,s)\,dy\,ds,
\end{equation}
Notice that $\widetilde F$ rescales exactly as~$F$. Moreover,
\begin{equation}
 \label{tFest}
|\widetilde F(x,t)|\le C_\eta|x|^{-\eta}t^{-(4-\eta)/2}, \qquad\hbox{for all $0\le \eta<\infty$}
\end{equation}
These are the same estimates as for $F$, but there is now no limitation
to the spatial decay rate (i.e., the restriction $\eta\le 4$ appearing in~\eqref{Fest}
can be removed).

Therefore, the space-time pointwise decay estimates for
$\widetilde B(\theta,u)(x,t)$ can be proved essentially in the same way as
in the previous Lemma.

The $L^1$-estimate (useful for estimating the $\Yy$-norm is straightforward:
$$ \bigl\|\widetilde B(\theta,u)(t)\bigr\|_1\le C\int_0^t (t-s)^{-1/2}\|u(s)\|_\infty\|\theta(s)\|_1\,ds
\le C\|u\|_{\mathcal{X}}\|\theta\|_{\mathcal{Y}}.$$
This allows us to conclude.

\endProof

%
%
%
%

\paragraph{\bf Proof of Proposition~\ref{prop1}.}
We need two elementary estimates on the linear heat equation. Namely,
\begin{equation}
 \|e^{t\Delta}\theta_0\|_{\Yy}\le C\Bigl(\|\theta_0\|_1+\hbox{ess}\sup_x|x|^{3}|\theta_0(x)|\Bigr)
\end{equation}
and
\begin{equation}
 \|e^{t\Delta}u_0\|_\Xx\le C\hbox{ess}\sup_{x} |x|\,|u_0(x)|.
\end{equation}
Both estimates immediately follow  from direct computations on the heat kernel
$g_t(x)=(4\pi\,t)^{-3/2}e^{-|x|^2/(4t)}$. 
(See, {\it e.g.\/} \cite{Bra04i, Miy00})
Here one only needs to use $|g_t(x)|\le C|x|^{-3}$ and the usual
$L^1$-$L^\infty$ estimates for $g_t$.
Letting $U=e^{t\Delta}u_0$ and $\Theta=e^{t\Delta}\theta_0$, by assumption~\eqref{smalla}
we get, for some $C>0$, $\|U\|_\Xx+\|\Theta\|_\Yy\le C\epsilon$.
The system~\eqref{IE} can be written in the abstract form~\eqref{ABS}.
By inequalities~\eqref{obm}, \eqref{bbe3} and \eqref{bte3}, all the assumptions
of Lemma~\ref{LKP} are satisfied provided $\epsilon>0$ is small enough.
The conclusion of Proposition~\ref{prop1} readily follows.

\endProof

\paragraph{\bf Proof of Part (a) of Proposition~\ref{prop2}.}
By construction, the solution $(u,\theta)$ of Proposition~\ref{prop1}
is obtained as the limit in $\Xx\times \Yy$ of the sequence $(u^k,\theta^k)$
defined in~\eqref{ites}.

By the first of assumptions~\eqref{adde}, and applying straightforward estimates on the heat kernel
(see also~\cite{Bra04i, Miy00}),
$|e^{t\Delta}u_0(x)|\le C(1+|x|)^{-a}$ and $|e^{t\Delta}u_0|\le C(1+t)^{-a/2}$.
(Here we need $0\le a<3$).
These two conditions imply in particular that $e^{t\Delta}u_0\in \Xx_a$.
Similarily  one deduces from the second inequality in~\eqref{adde}
that $e^{t\Delta}\theta_0\in \Yy_b$
(when $b=3$ one here needs also $\theta_0\in L^1$).

By estimate~\eqref{bte1}
and assumption~\eqref{smalla}, for all $k=1,2\ldots$
we get 
$$\|\theta^{k+1}\|_{\Yy_b}\le \|e^{t\Delta}\theta_0\|_{\Yy_b}
+C\epsilon\|\theta^k\|_{\Yy_b}.$$
If $\epsilon>0$ is small enough then $C\epsilon<1$
(the size of the admissible $\epsilon$ thus depend on $b$ and, as we will see later,
also on~$a$).
Iterating this inequality shows that the sequence $(\theta^k)$
is bounded in $\Yy_b$.

Combining estimates~\eqref{ltx} with~\eqref{bbe1}
we get
$$
\|u^{k+1}\|_{\Xx_a}\le \|e^{t\Delta}u_0\|_{\Xx_a}+C\epsilon\|u^k\|_{\Xx_a}+C\|\theta_k\|_{\Yy_3}.
$$
Assuming $C\epsilon<1$, we deduce from the boundness of $(\theta^k)$
in $\Yy_3$ that $(u^k)$ is bounded in~$\Xx_a$.
Thus, the solution $(u,\theta)$ belongs to $\Xx_a\times \Yy_b$
and the first part of Proposition~\ref{prop2} follows.

\bigskip
\paragraph{\bf Proof of Part (b) of Proposition~\ref{prop2}.}
The proof of the second part of Proposition~\eqref{prop2}
is quite similar but relies on the use of slightly different function spaces.
So, let $a\ge2$. We define $\tXx_a$ as the Banach space of divergence vector fields $u=u(x,t)$
such that, for some $C>0$,
\begin{equation}
\label{deca1bis}
|u(x,t)|\le C\inf_{0\le \eta\le a} |x|^{-\eta}(1+t)^{(\eta-2)/2}.
\end{equation}
For $b\ge 4$ we define the space $\tYy_b$ of functons $\theta\in L^\infty_t(L^1_1)$ satisfying the estimates
\begin{equation}
\label{edeca11bis}
|\theta(x,t)|\le C\inf_{0\le \eta\le b} |x|^{-\eta}(1+t)^{(\eta-4)/2}.
\end{equation}
Such spaces are equipped with their natural norms.

Notice that the spaces $\tXx_a$ and $\tYy_b$ differ from the their counterparts $\Xx_a$ and $\Yy_b$
only by the fact that the time decay conditions are slightly more stringent in the former case.

The counterpart of estimates~\eqref{elx} are
\begin{equation}
 \label{celx}
\begin{split}
&\|u(t)\|_{L^{p,q}}\le C\|u\|_{\tXx_a}
 (1+t)^{\frac{1}{2}(\frac{3}{p}-2)}, \qquad \max\{1,\textstyle\frac{3}{a}\}<p\le\infty,\\
&\|\theta(t)\|_{L^{p,q}}\le C\|\theta\|_{\tYy_4} (1+t)^{\frac{1}{2}(\frac{3}{p}-4)},
\qquad 1\le p\le\infty.\\
\end{split}
\end{equation}
The first of~\eqref{celx} can be completed by
\begin{equation}
 \label{celx2}
\|u(t)\|_1\le C\|u\|_{\widetilde{\mathcal{X}}_a}(1+t)^{-1+\frac{3}{a}},
\qquad\hbox{if $3<a<4$.}
\end{equation}
This last estimate follows immediately by splitting the integral $\int |u|$
into the regions $|x|\ge t^{1/a}$ and $|x|\le t^{1/a}$.

We also notice the continuous embedding
\begin{equation}
\label{embp}
\widetilde\Yy_b \subset L^\infty_t(L^1_1), \qquad \hbox{if $b>4$}
\end{equation}
that follows easily by splitting the integral 
$\int (1+|x|)|\theta(x,t)|\,dx$ into $\int_{|x|\le (1+t)^{1/2}}\dots$ and $\int_{|x|\ge (1+t)^{1/2}}\dots$,
and using the bound $|\theta(x,t)|\le C(1+t)^{-2}$ for the first term and $|\theta(x,t)|\le C|x|^{-b}$
for the second one.

\begin{lemma}
 \label{adl1}
Let $2\le a<4$. If $\int\theta(t)\,dx=0$ for all $t$, then for some $C>0$,
\begin{equation}
 \label{ae1}
\|L(\theta)\|_{\tXx_a}\le C\|\theta\|_{\tYy_4}.
\end{equation}
\end{lemma}
Assume, without restriction, $\|\theta\|_{\tYy_4}=1$.
The time decay estimate 
$$\|L(\theta)\|_\infty=\Bigl\|\int_0^t\K(t-s)*\theta(s)\,ds\Bigr\|_\infty
\le C(1+t)^{-1}$$
immediately follows by the $L^2$-$L^2$ Young inequality.
It only remains to prove that $L(\theta)$ can be bounded by $C|x|^{-a}(1+t)^{(a-2)/2}$
for all $(x,t)$ belonging to the parabolic region $|x|\ge 2\sqrt t$.
Thus, we decompose
\begin{equation}
\label{aee2}
L(\theta)=I_1+I_2+I_3=(I_{1,1}+I_{1,2}+I_{1,3})+I_2+I_3.
\end{equation}
where $I_1$, $I_2$ and $I_3$ are as in Lemma~\ref{L34}
and and the terms contributing to $I_1$
are defined below.
First,
\begin{equation*}
 \bigl(I_{1,1}\bigr)_j(x,t)\equiv\int_0^t\K_{j,3}(x,t-s)\int\theta(y,s)\,dy\,ds=0
\end{equation*}
by the zero-mean assumption on~$\theta$.
Next
\begin{equation*}
 \bigl(I_{1,2}\bigr)_j(x,t)\equiv-\int_0^t\K_{j,3}(x,t-s)\int_{|y|\ge |x|/2}\theta(y,s)\,dy\,ds
\end{equation*}
and
\begin{equation*}
 \bigl(I_{1,3})_j(x,t)\equiv-\int_0^t\!\!\int_0^1\biggl(\int_{|y|\le |x|/2} \nabla \K_{j,3}
(x-\lambda y,t-s)\,dy\biggr)
 \cdot y\,\theta(y,s)\,d\lambda\,ds,
\end{equation*}
where we have used the Taylor formula to write the difference $\K_{j,3}(x-y,t-s)-\K_{j,3}(x,t-s)$.

From the bounds $|K(x,t)|\le C|x|^{-3}$ and $|\nabla K(x,t)|\le C|x|^{-4}$
we get
\begin{equation*}
 |I_{1,2}|+|I_{1,3}|\le C|x|^{-4}t,
\end{equation*}
where we used the continuous embedding of $L^1_1$ into $\tYy_4$
that follows from the definition of such space.
The above pointwise estimate is even better, in our parabolic region $|x|\ge 2\sqrt t$,
than what we actually need.

Next,
 $$(I_2)_j=\int_0^t\!\!\int_{|x-y|\le |x|/2}\K_{j,3}(x-y,t-s)\theta(y,s)\,dy\,ds$$
can be bounded as follows
\begin{equation*} 
\begin{split}
|I_2| &\le 
C|x|^{-4}\int_0^t\!\!\int_{|x-y|\le |x|/2}|\K_{j,3}(x-y,t-s)|\,dy\,ds\\
& \le C|x|^{-4}t\log(|x|/\sqrt t)\\
&  \le  C|x|^{-a}t^{(a-2)/2}
\end{split}
\end{equation*}
for all $(x,t)$ such that $|x|\ge 2\sqrt t$ (recall that $2\le a<4$).

Moreover, using again $|K(x-y,t-s)|\le C|x-y|^{-3}$, shows that
$$(I_3)_j=\int_0^t\!\!\int_{|y|\ge |x|/2,\; |x-y|\ge |x|/2}\K_{j,3}(x-y,t-s)\theta(y,s)$$
can be bounded by 
\begin{equation*}
 C|x|^{-4}t.
\end{equation*}

So far, we proved that
$$ |L(\theta)|(x,t) \le C|x|^{-a}t^{(a-2)/2}\le C|x|^{-a}(1+t)^{(a-2)/2}$$
for all $(x,t)$ such that $|x|\ge 2\sqrt t$.
Our previous $L^\infty$-bound on $L(\theta)$ implies the validity
of such estimate in the region $|x|\le 2\sqrt t$. We thus conclude
that $L(\theta)\in \tXx_a$
and Lemma~\ref{adl1} follows.

\endProof

Next Lemma is a simple variant of Lemma~\ref{buv}.

\begin{lemma}
 \label{borl}
Let $2\le a<4$ and $b\ge 4$. Then, for some constant $C>0$,
\begin{equation}
\label{bbet} 
\bigl\| B(u,v)\bigr\|_{\tXx_a}\le C\bigl\|u\bigr\|_{\Xx}\bigl\|v\bigr\|_{\tXx_a}\,\,
\end{equation}
\end{lemma}
and 
\begin{equation}
\label{bbetbis} 
\bigl\| \widetilde B(u,\theta)\bigr\|_{\tYy_b}\le C\bigl\|u\bigr\|_{\Xx}\bigl\|\theta\bigr\|_{\tYy_b}\,\,
\end{equation}

\Proof
This Lemma can be easily proved following the
steps of estimates~\eqref{bbe1} and~\eqref{bte1}.
We thus skip the details.

\endProof

The last estimates that we need, concern the heat equation.
The computations are straightforward (see \cite{Bra04i, Miy00}).
Recall that $2\le a<4$ and $a\not=3$. Moreover, we assumed
$|u_0(x)|\le C(1+|x|)^{-a}$.
When $2\le a<3$ then, as we already observed, 
$|e^{t\Delta}u_0(x)|\le C(1+|x|)^{-a}$ and $|e^{t\Delta}u_0|\le C(1+t)^{-a/2}$.
In fact this estimate remains valid also for $3<a<4$ (here one uses that
$u_0$ is integrable and divergence free, and so $\int u_0=0$).
Thus, in particular, $e^{t\Delta}u_0\in \tXx_a$.

In the same way, one proves that by our assumptions
$e^{t\Delta}\theta_0\in \tYy_b$ for $b\ge4$.

Therefore, going back to the approximation scheme~\eqref{ites}
and arguing as in the proof of Part (a) of Proposition~\eqref{prop2}
we see that the sequence $(\theta^k)$ is bounded in $\tYy_b$
and $(u^k)$ is bounded in $\tXx_a$.
Part (b) of Proposition~\eqref{prop2}  follows.

\endProof

\section{Asymptotic profiles and decay of strong solutions}
\label{sec:eight}

We denote by $E(x)=\frac{c}{|x|}$ the fundamental solution of the Laplacian in~$\R^3$
and by $\bigl(E_{x_j,x_k}\bigr)(x)$ its second order derivatives for $x\not=0$.
Notice that $E_{x_x,x_3}$ is a homogeneous function of degree~$-3$.
Next lemma describes the asymptotic profile for $L(\theta)(x,t)$ as $|x|\to\infty$,
by establishing that
$$L(\theta)(x,t)\simeq  \biggl(\int\theta_0\biggr)\,t\,\bigl(E_{x_j,x_3}\bigr)(x),
\qquad\hbox{as $|x|>\!\!\!>\sqrt t$}.$$

\begin{lemma}
 \label{prol}
Let $\theta=\theta(x,t)$ be any function satisfying the pointwise estimates~\eqref{deca11},
for some~$3<b<4$ and such that  $\int\theta(t)=\int\theta_0$ for all $t\ge0$.
Then the $j$-component $L(\theta)_j(t)$ of $L(\theta)$
can be decomposed as
\begin{equation}
 \label{proe}
L(\theta)_j(x,t)=\biggl(\int\theta_0\biggr)\,t\,\bigl(E_{x_j,x_3}\bigr)(x)\,+\,
\mathcal{R}'_j(x,t)
\qquad(j=1,2,3),
\end{equation}
where the remainder function $\mathcal{R'}$ satisfies,
\begin{equation}
|\mathcal{R'}(x,t)|\le C|x|^{-b}\,t^{(b-1)/2}\,\log(|x|/\sqrt t),
\qquad\forall\, (x,t) \;\hbox{ s.t. }\; |x|\ge 2\sqrt t.
\end{equation}

In particular,  in the region $|x|>\!\!\!>\sqrt t$,  one has
$$|\mathcal{R'}(x,t)|<\!\!\!< Ct\,|E_{x_j,x_k}(x)|$$
along almost all directions.
\end{lemma}

\begin{remark}
This idea of obtaining informations on the large time behavior 
of solutions by first studying their behavior in the
parabolic region $|x|>\!\!\!>\sqrt t$ comes from~\cite{BraV07}.
\end{remark}

\medskip
\Proof
We go back to the decomposition~\eqref{aee2}
of $L(\theta)$, as done in Lemma~\ref{adl1}.
We now treat $I_2$ using  the estimate
$|\theta(x,t)|\le C|x|^{-b}\,(1+t)^{(b-3)/2}$.
This yields to the inequality, valid for $|x|\ge 2\sqrt t$, 
\begin{equation*}
 |I_2|(x,t)\le C|x|^{-b}(1+t)^{(b-1)/2}\log(|x|/\sqrt t).
\end{equation*}
With the bound on $\theta$ we obtain also
\begin{equation*}
 |I_3|(x,t)\le C|x|^{-b}(1+t)^{(b-1)/2}.
\end{equation*}
Now recall that $I_1=I_{1,1}+I_{1,2}+I_{1,3}$, where
\begin{equation}
 \bigl(I_{1,3})_j(x,t)=-\int_0^t\!\!\int_0^1\biggl(\int_{|y|\le |x|/2} \nabla \K_{j,3}(x-\lambda
 y,t-s)\,d\lambda\biggr)
 \cdot y\theta(y,s)\,dy\,ds.
\end{equation}
and that $\nabla \K$ satisfies the estimate $|\nabla \K(x,t)|\le C|x|^{-4}$.
Then, since $3<b<4$,
\begin{equation*}
\begin{split}
|I_{1,3}|(x,t) 
&\le C|x|^{-4}\int_0^t\!\!\int_{|y|\le |x|/2}|y|\,\,|\theta(y,s)|\,dy\,ds\\
&\le C|x|^{-b}\,(1+t)^{(b-1)/2}.
\end{split}
\end{equation*}

The estimate for $(I_{1,2})_j=-\int_0^t\K_{j,3}(x,t-s)\int_{|y|\ge |x|/2}\theta(y,s)\,dy\,ds,$ is straightforward:
\begin{equation*}
|I_{1,2}|(x,t) \le C|x|^{-b}\,(1+t)^{(b-1)/2}.
\end{equation*}

Finally,
since the mean of~$\theta$ remains constant in time,
\begin{equation*}
 \bigl(I_{1,1}\bigr)_j(x,t)=\int_0^t\K_{j,3}(x,t-s)\int\theta(y,s)\,dy\,ds
=\biggl(\int \theta_0\biggr)\int\K_{j,3}(x,t-s)\,ds.
\end{equation*}

But the following decomposition of the kernel~$\K$, established in \cite{Bra08},  holds :
\begin{equation*}
\K_{j,k}(x,t)= E_{x_j,x_k}(x)+|x|^{-3}\Psi_{j,k}(x/\sqrt t), \qquad j,k=1,2,3
\end{equation*}
where $\Psi_{j,k}$ is fast decaying:  $|\Psi(y)|\le Ce^{-c|y|^2}$ for all $y\in\R^3$ and some constats $c,C>0$.
Hence, we can estimate $|\Psi(y)|\le C|y|^{-b+3}$.

We now define~$\mathcal{R}'(x,t)$ through the relation
$$ L(\theta)(x,t)=  \biggl(\int\theta_0\biggr)\,t\,\bigl(E_{x_j,x_3}\bigr)(x)+\mathcal{R}'(x,t)$$
and all the previous estimates imply
$|\mathcal{R}'(x,t)|\le C|x|^{-b}(1+t)^{(b-1)/2}\log(|x|/\sqrt t)$
for all $|x|\ge 2\sqrt t$.

\endProof

In the case $\int \theta_0=0$, we can use  the following variant
of Lemma~\ref{prol}.

\begin{lemma}
 \label{prol2}
Let $\theta=\theta(x,t)$ be any function satisfying the second of~\eqref{deca2},
for some~$4<b<5$ and such that  $\int\theta(t)=\int\theta_0=0$ for all $t\ge0$.
Then the $j$-component $L(\theta)_j(t)$ of $L(\theta)$
can be decomposed as
\begin{equation}
 \label{proebis}
L(\theta)_j(x,t)=
-\nabla E_{x_jx_3}(x)\cdot\biggl(\int_0^t\!\!\int y\,\theta(y,s)\,dy\,ds\biggr)+{\mathcal{R}''}_{j}(x,t)
\end{equation}
with 
\begin{equation} 
\label{wtr}
|{\mathcal{R}''}|(x,t)\le C|x|^{-b}(1+t)^{(b-2)/2}\log(|x|/\sqrt t)
\end{equation}
for all $(x,t)$ such that $|x|\ge 2\sqrt t$.
%
\end{lemma}

\medskip
\Proof
We only have to reproduce the proof of the previous Lemma with
slight modification.
Using the estimate
$|\theta(x,t)|\le C|x|^{-b}\,(1+t)^{(b-4)/2}$ we now obtain,
for $|x|\ge 2\sqrt t$, 
\begin{equation*}
 |I_2|(x,t)\le C|x|^{-b}(1+t)^{(b-2)/2}\log(|x|/\sqrt t),
\end{equation*}
next
\begin{equation*}
 |I_3|(x,t)\le C|x|^{-b}(1+t)^{(b-2)/2}.
\end{equation*}
and 
\begin{equation*}
 |I_{1,2}|(x,t)\le C|x|^{-b}(1+t)^{(b-2)/2}.
\end{equation*}

Next, by the vanishing mean condition $I_{1,1}=0$.
It remains to treat $I_{1,3}$.
We can decompose $I_{1,3}$, whose $j$-component we recall is
$$\int_0^t\!\!\int_{|y|\le |x|/2} [\K_{j,3}(x-y,t-s)-\K_{j,3}(x,t-s)]\theta(y,s)\,dy\,ds,$$
 into the sum of three more terms
$$I_{1,3}=I_{1,3,1}+I_{1,3,2}+I_{1,3,3}.$$
Such decomposition is performed exactly in the way we did in the proof of Lemma~\ref{prol}.
The $j$-component of the first term is thus
\begin{equation}
\label{qp6} 
(I_{1,3,1})_j(x,t)=-\int _0^t \nabla\K_{j,3}(x,t-s)\cdot\int y\,\theta(y,s)\,dy\,ds.
\end{equation}
The second term,
$$(I_{1,3,2})_j(x,t)=\int_0^t \nabla\K_{j,3}(x,t-s)\cdot\int_{|y|\ge |x|/2} y\theta(y,s)\,dy\,ds,$$
 can be bounded by the right-hand side of~\eqref{wtr}
using $|\nabla \K(x,t)|\le C|x|^{-4}$ and $|\theta(x,t)|\le C|x|^{-b}(1+t)^{(b-4)/2}$.
Next, the $j$-component of $I_{1,3,3}$,
$$-\int_0^t\!\!\int_{|y|\le |x|/2} [\K_{j,3}(x-y,t-s)-\K_{j,3}(x,t-s)+\nabla\K(x,t-s)\cdot y]\theta(y,s)\,dy\,ds,$$
can be treated with the Taylor formula.
The simple estimate $|\nabla^2_x\K(x,t)|\le C|x|^{-5}$
allows us to see that also $I_{1,3,3}$ is bounded by the right-hand side of~\eqref{wtr}.
Therefore, both $I_{1,3,2}$ and $I_{1,3,3}$ can be included into the remainder term ${\mathcal{R}''}(x,t)$.

Let us go back to~\eqref{qp6}.
As shown in~\cite{Bra08}, the following decomposition holds true:
\begin{equation}
\nabla\K_{j,k}(x,t)=\nabla E_{x_jx_k}(x)+|x|^{-4}\widetilde\Psi(x/\sqrt t),
\end{equation}
with $|\widetilde\Psi(y)|\le Ce^{-c|y|^2}$ for some constants $C,c>0$ and all $y\in\R^3$.
In particular we can estimate $|\widetilde\Psi(y)|\le C|y|^{-(b-4)}$.
On the other hand, $\int |y|\,|\theta(y,s)|\,dy$ is uniformly bounded because of the 
embedding~\eqref{embp}.
This shows that  $(I_{1,3,1})_j(x,t)$ can be written in the region $|x|\ge 2\sqrt t$
as in the right-hand side of~\eqref{wtr} (even without logarithmig factors).

This finally gives~\eqref{proebis}.

\endProof

\bigskip
We can now establish  our main results as simple corollaries:

\bigskip

\paragraph{\bf Proof of Theorem~\ref{theor1}, part (a).}
Let $(u,\theta)$ be a mild solution of the system~\eqref{IE}, satisfying the pointwise decay
 estimates~\eqref{deca1}-\eqref{deca11}, with~$a>\frac{3}{2}$ and $b>3$.
Recall that the spaces $\Xx_a$ and $\Yy_b$ decrease as $a$ and $b$ grow.
Without restriction we can then assume $\frac32<a<3$ and $3<b<4$ in our  calculations.
According to our notations, we can write
\begin{equation}
\label{ueqa}
u(x,t)=e^{t\Delta}u_0(x)+B(u,u)(x,t)+ L(\theta)(x,t).
\end{equation}
By estimate~\eqref{bbe2}, owing to the condition $a>\frac{3}{2}$,
we have
\begin{equation*}
\lim_{\frac{|x|}{\sqrt t}\to\infty} \frac{B(u,u)(x,t)}{t|x|^{-3}}=0.
\end{equation*}
Therefore $B(u,u)$ can be included inside the remainder term in the asymptotic profile
of~$u$ for $\frac{|x|}{\sqrt t}\to\infty$.
Moreover Lemma~\eqref{prol} and the condition $b>3$ guarantee that
\begin{equation*}
L(\theta)(x,t)=\Bigl(\int \theta_0\Bigr)\,t\,\nabla E_{x_3}(x) + o\bigl(t|x|^{-3}\bigr),
\qquad\hbox{as $\frac{|x|}{\sqrt t}\to\infty$}.
\end{equation*}
This yields the asymptotic profile~\eqref{coroc} for~$u$.

Let us prove here also the claim in Part~(a) of Remark~\ref{rem-ullar}.
As usual, we denote $g_t(x)=(4\pi t)^{-3/2}e^{-|x|^2/(4t)}$ the standard gaussian.
Under the additional assumption $|u_0|\le C|x|^{-3}$,
we have $\int_{|y|\ge |x|/2} g_t(x-y)|u_0(y)|\,dy\le C|x|^{-3}$.
Moreover, by the third of~\eqref{smalla}, we have also, e.g.,  $|u_0(x)|\le C|x|^{-2}$.
So 
$\int_{|y|\le |x|/2} g_t(x-y)|u_0(y)|\,dy\le C|x|\sup_{|z-x|\le |x|/2} g_t(z)\le \sqrt t\,|x|^{-3}$.
Combining these estimates we get
$$ |e^{t\Delta}u_0(x)| <\!\!\!< t|x|^{-3}, \qquad\hbox{as $t>\!\!\!>1$}.$$
Thus, when $\int\theta_0\not=0$, under the additional condition,
$|u_0(x)|\le C|x|^{-3}$, for $|x|>\!\!\!>\sqrt t>\!\!\!>1$, the solution $u$ behaves like
$\bigl(\int u_0\bigr)t\nabla E_{x_3}(x)$ along almost all directions.
More precisely, using spherical coordinates and letting 
$x=\rho\omega$, with $\rho>0$, for almost all $\omega$ in the unit sphere, we have
$$\lim _{\rho,t\to\infty} {u_j(x,t)}/\Bigl({(\textstyle\int u_0)t E_{x_j,x_k}(x)}\Bigr)=1.$$

\bigskip

\paragraph{\bf Proof of Theorem~\ref{theor1}, part (b).}
We can assume without restictrions for the calculations below that that 
$2<a<3$, and $b>4$. Moreover, by our assumption $\int\theta_0=0$.
The velocity field $u$ belongs to $\tXx_2$.
Then 
$$|B(u,u)|(x,t)\le C|x|^{-4}t^{1/2}$$
as it can be proved easily by
splitting $B(u,u)$ into $\int_0^t\!\!\int_{|y|\le |x|/2}$ and $\int_0^t\!\!\int_{|y|\ge |x|/2}$
and using $|u|\le C|x|^{-2}$ and $\|u(t)\|\le Ct^{-1/4}$ (see~\eqref{celx}).
We get easily a bound for both terms in the region $|x|\ge 2\sqrt t$, implying
\begin{equation*}
 \lim_{t,\,\frac{|x|}{\sqrt t}\to\infty} \frac{|B(u,u)|(x,t)}{t|x|^{-4}}=0
\end{equation*}
Thus, $B(u,u)$ can be included inside the remainder term.
Applying now Lemma~\ref{prol2} yields the asymptotic expansion~\eqref{putz}.

Under the additional assumption $|u_0(x)|\le C|x|^{-4}$,
we get, by~\eqref{smalla}, $|u_0(x)|\le C|x|^{-5/2}$.
If we use the bound $|g_t(x)|\le Ct^{3/4}|x|^{-9/2}$ and the usual $L^1$-estimate
for~$g_t$ we get $|e^{t\Delta}a(x)|\le C(1+t^{3/4})|x|^{-4}$.
Thus,
$$ |e^{t\Delta}u_0(x)| <\!\!\!< t|x|^{-4}, \qquad\hbox{as $t>\!\!\!>1$}$$
and the last claim (made rigorous exactly as above)
of the theorem follows.

\endProof

\medskip
We now deduce from Theorem~\ref{theor1}
sharp upper and lower bound estimates in $L^p$-spaces.

\medskip

\paragraph{\bf Proof of Corollary~\ref{coro1}, part (a).}
The upper bounds are simple:
indeed, appying the arguments that we used in the proof of Lemma~\ref{emb}
to $(1+|\!\cdot\!|)^r|u(\cdot,t)|$ instead of~$u$ (and putting $p=q$, in a such way
that Lorentz spaces boil down to the usual Lebesgue spaces) gives the result.

We now discuss lower bounds.
By the proof of Lemma~\ref{prol} and of Theorem~\ref{theor1},
we can find an exponent $\eta>0$ (any $0<\eta<\min\{2a-3,b-3,1\}$ will do)
such that, for $j=1,2,3$,
\begin{equation*}
 \begin{split}
  &|u_j(x,t)-e^{t\Delta}(u_{0,j})(x)|\\
&\qquad \ge t\,\Bigl|\textstyle\int\theta_0\Bigr| \,\, |E_{x_j,x_3}(x)|
 -Ct|x|^{-3}\bigl(\frac{|x|}{\sqrt t}\bigr)^{-\eta}
 \end{split}
\end{equation*}
provided $|x|\ge A\sqrt t$, and $A>0$ is large enough.

Consider the parabolic region $\mathcal{D}_{A,t}=\{(x,t)\colon |x|\ge A\sqrt t\}$.
For $1<p<\infty$, we denote by
$\|\cdot\|_{L^p_r(\mathcal{D}_{A,t})}$ the norm
$$ \|f\|_{L^p_r(\mathcal{D}_{A,t})}=\biggl(\int_{\mathcal{D}_{A,t}} |f(x)|^p(1+|x|)^{rp}\,dx\biggr)^{1/p}.$$
Then for all~$t\ge1$ and $1<p<\infty$, $r\ge0$ such that $r+\frac{3}{p}<3$,
\begin{equation*}
 \begin{split}
  &\Bigl\| u_j(t)-e^{t\Delta}u_{0,j}\Bigr\|_{L^p_r(\mathcal{D}_{A,t})}\\
  &\qquad \ge 
  C\,t\,\Bigl|\textstyle\int\theta_0\Bigr|\,\,
 \Bigl\|E_{x_j,x_3}\Bigr\|_{L^p_r(\mathcal{D}_{A,t})}
  -C\,t \Bigl\| |\!\cdot\!|^{-3} \bigl(\frac{|\cdot|}{\sqrt t}\bigr)^{-\eta}
                    \Bigr\|_{L^p_r(\mathcal{D}_{A,t})}\\
&\qquad\ge Ct^{\frac{1}{2}(r+\frac{3}{p}-1)}  A^{r+\frac{3}{p}-3}
\biggl(C'\Bigl|\textstyle\int \theta_0\Bigr|-A^{-\eta}\Biggr). 
 \end{split} 
\end{equation*}
This shows that it is possible to define a continuous function $\phi\colon[0,\epsilon^2]\to\R^+$,
(where $\epsilon>0$ is the constant of Proposition~\ref{prop1})
such that $\phi(0)=0$, $\phi$ is strictly positive outside the origin, and satisfying
\begin{equation}
\label{dah}
\begin{split} 
 \Bigl\|u(t)-e^{t\Delta}u_0\Bigr\|_{L^p_r} 
 &\ge \Bigl\|u(t)-e^{t\Delta}u_0\Bigr\|_{L^p_r(\mathcal{D}_{A,t})}\\
 &\ge \phi\bigl(|\textstyle\int\theta_0|\bigr)\,t^{\frac{1}{2}(r+\frac{3}{p}-1)} 
\end{split}
\end{equation}
for all $t>0$ large enough.
By comparing the two terms inside the parentheses in the inequality above,
{\it i.e.\/}, by taking $A$ such that $A^{-\eta}\le \frac{C'}{2}|\int\theta_0|$, 
we get an explicit behavior for~$\phi(\sigma)$ near zero, namely, 
$\phi(\sigma)\sim c\sigma^{1+\frac{1}{\eta}(3-r-\frac{3}{p})}$, as $\sigma\to0^+$.
with $c>0$ small enough.

We now restrict us to the smaller range $0\le r+\frac3p<\min(3,a)$,
always with $r\ge 0$ and $1<p<\infty$. Let us compute the $L^p_r$-norm
of~$e^{t\Delta}u_0$.
From $|u_0(x)|\le C\min\{|x|^{-1},|x|^{-a}\}$ we obtain,
for $t\ge1$, $e^{t\Delta/2}|u_0|(x)(1+|x|)^r\le C(1+|x|)^{-(a-r)}$.
Applying the semigroup property of the heat kernel, we get, for $t\ge2$,
$e^{t\Delta}|u_0|(x)(1+|x|)^r\le e^{t\Delta/2}(1+|x|)^{-(a-r)}$.
Computing the $L^p$-norm of this quantity, we deduce
\begin{equation}
\label{dahli}
 \begin{split}
 \|e^{t\Delta}u_0\|_{L^p_r}&\le \|g_{t/2}*(1+|\cdot|)^{-(a-r)}\|_p\\
&\le C\|g_{t/2}\|_{L^{\alpha,p}}\\
 &=Ct^{-\frac{1}{2}(a-r-\frac3p)}.
\end{split}
 \qquad \hbox{with $1+\frac{1}{p}=\frac{1}{\alpha}+\frac{a-r}{3}$}
\end{equation}
In this computation, $L^p_r$ denotes as usual the weighted $L^p$ space, whereas $L^{\alpha,p}$ is a Lorentz space.
Here we made use of Young convolution inequality, generalised to Lorentz spaces
(see~\cite[Prop. 2.4]{Lem02}).

By comparing the large time behavior of the RHS in expressions~\eqref{dah}-\eqref{dahli},
we deduce the lower bound
\begin{equation*}
 \|u(t)\|_{L^p_r}
\ge \textstyle\frac{1}{2}\phi\bigl(|\textstyle\int\theta_0|\bigr)
\,t^{\frac{1}{2}(r+\frac{3}{p}-1)},
\qquad \hbox{for all $t\ge t_0$}
\end{equation*}
where $t_0>0$ is some constant 
depending on all the parameters and the initial data, but independent
on~$t$.

\endProof

\medskip

\paragraph{\bf Proof of Corollary~\ref{coro1}, Part (b).}
The estimate from above follows applying inequalities~\eqref{celx} and~\eqref{celx2}
to $(1+|\cdot|)^r u$.

Let us now estimate $\|u(t)\|_{L^p_r}$ from below.
The proof is based on the asymptotic expansion~\eqref{putz}.
Computing the third order derivatives outside the origin of the fundamental solution $E$ of
$-\Delta$ in $\R^3$, i.e., $E(x)=\frac{C}{|x|}$,
shows that (see also~\cite[Eq. (9b)]{Bra08})
\begin{equation}
 \label{fracf}
E_{x_j,x_h,x_k}(x)=\frac{\Gamma\Bigl(\frac52\Bigr)}{\pi^{3/2}}\cdot
\frac{\sigma_{j,h,k}(x)|x|^2-5x_jx_hx_k}{|x|^7},
\end{equation}
with $\sigma_{j,h,k}(x)=\delta_{j,h}x_k+\delta_{h,k}x_j+\delta_{k,j}x_h$.
It is now easy to see that  the expression 
$\sum_{h,k=1}^d E_{x_j,x_h,x_k}(x) M_{h,k}$
identically vanishes
if and only if $M_{h,k}$ is a scalar multiple of the identity matrix.
Let ${\bf m}(t)=\int_0^t\!\!\int y\theta(y,s)\,dy\,ds$.
We deduce that
the homogeneous function of degree $-4$
$$ \nabla E_{x_j,x_3}(x)\cdot {\bf m}(t)$$
is identically zero, for any fixed~$t$, if and only if $|{\bf m}(t)|=0$.

If $\widetilde{\boldsymbol{m}}=\liminf_{t\to\infty} \left| \frac{1}{t}\int_0^t\!\!\int y\theta(y,s)\,dy\,ds\right|\not=0$,
then there exists $c>0$ such that, for all $t$ sufficiently large,
$$ \|\nabla E_{x_j,x_3}\cdot{\bf m}(t)\|_{L^p_r(\mathcal{D}_{A,t})}\ge 
c\,\widetilde{\boldsymbol{m}}\,t^{\frac{1}{2}(r+\frac{3}{p}-2)}.$$

The condition on the remainder~$\widetilde{\mathcal{R}}$ obtained in
Theorem~\ref{theor1} then implies, some constant $c>0$ and for all
$t$ large enough,
\begin{equation*}
 \|u(t)-e^{t\Delta}u_0\|_{L^p_r} \ge ct^{\frac12(r+\frac3p-2)}
\end{equation*}
with $r\ge0$, $1<p<\infty$ tels que $r+\frac3p<4$.
It remains to prove that, when $r+\frac3p<\min\{a,4\}$ then we have
\begin{equation}
\label{thee}
\|e^{t\Delta}u_0\|_{L^p_r} = o(t^{\frac12(r+\frac3p-2)}), \qquad
\hbox{for $t\to\infty$}.
\end{equation}
Recall that we assumed $a>2$. When $2<a<3$,  we can simply use
inequality~\eqref{dahli}.
When $a=3$ there is nothing to prove because we reduce to the previous case
by picking $a'$, with $r+\frac3p<a'<3$.
So, consider now $3<a<4$.
In this case $u_0$ is integrable and $\int u_0=0$ by the divergence-free condition.
Then, for $t\ge1$, $e^{t\Delta/2}|u_0|(x)(1+|x|)^r\le C(1+|x|)^{-(a-r)}$.
Thus,
\begin{equation}
\label{dahli2}
 \begin{split}
 \|e^{t\Delta}u_0\|_{L^p_r}&\le \|g_{t/2}*(1+|\cdot|)^{-(a-r)}\|_p
\end{split}
\end{equation}
When $a-r<3$ we can apply Young inequality in the same way as before
and still obtain estimate~\eqref{dahli}.
When $a-r>3$, the above quantity is bounded by
$Ct^{-\frac12(\frac3p-3)}$ (and by $C_\eta t^{-\frac12(\frac3p-3+\eta)}$
for all $\eta>0$ when $a-3=3$).
In any case,~\eqref{thee} holds true.
This establishes estimates~\eqref{sdabz2}

\endProof

\medskip
Theorem~\ref{theor1} has another interesting consequence, that clarifies the importance
of the restriction $r+3/p<3$ in our previous statements.

\begin{corollary}
 \label{reac}
Let $(\theta,u)$ be a solution as in Part~(a) of Theorem~\ref{theor1}.
We assume, in addition, that $\int \theta_0\not=0$ and 
that the initial velocity satisfies $|u_0(x)|\le C(1+|x|)^{-a}$,
for some $a>3$.
Then for all
$$ r\ge 0, \qquad 1\le p<\infty,\qquad r+\frac{3}{p}\ge 3$$
and for all $t>0$ we have
\begin{equation}
 \|u(t)\|_{L^p_r}=\infty.
\end{equation}
\end{corollary}

\Proof
Indeed, $E_{x_j,x_3}(x)$ is a homogeneous function of degree~$3$, smooth outside the origin.
Then we can find an open conic set~$\Gamma\subset \R^3$ such that
$|E_{x_j,x_3}(x)|\ge C|x|^{-3}>0$ for all $x\in \Gamma$, $x\not=0$, for some $C>0$.
Indeed, by Eq.~\eqref{coroc},
for $|x|\ge A\sqrt t$ with $A>0$ large enough and $x\in \Gamma$, and $t>0$,
we have
$$|u_j(x,t)|\ge \frac{C}{2}\,t\,\biggl|\int\theta_0\biggr| |x|^{-3}.$$
This implies that $u_j$ has an infinite $\|\cdot\|_{L^p_r(\mathcal{D}_A\cap\Gamma)}$\,-norm.

\endProof

\section{Additional remarks and comments}
\label{sec:ee8}

In this section we collect a few technical remarks on the main results.
These are essentially small variants of  our statements that can be easily proved
with minor modifications to the proofs.

\begin{remark} [on Theorem~\ref{th:belo}]
 Part (b) of Theorem~\ref{th:belo} can be streghened as follows.
Under the same assumptions on $\theta_0$
and replacing the assumption  $u_0\in L^2_{\boldsymbol{\sigma}}\cap L^{3/2}$ with the weaker
condition
$\|e^{t\Delta}u_0\|^2\le C(1+t)^{-s}$,
for some $s\ge0$, 
(when $u_0\in L^2_{\boldsymbol{\sigma}}\cap L^{3/2}$, this holds with $s=1/2$),
we have
\begin{equation}
 \|u(t)\| \le C(1+t)^{-s^*}, \qquad\hbox{with $s^*=\min(s,1/2)$}.
\end{equation}
The above decay condition on $e^{t\Delta}u_0$ could also be restated in terms of Besov spaces.
\end{remark}

\begin{remark} [on Theorem~\ref{theor1}]
The second term in the RHS of~\eqref{putz} is bounded by $C|x|^{-4}t$.
If one is interested in studying the asymptotic behavior of~$u$
only as $\frac{|x|}{\sqrt t}\to\infty$
(with $t>0$ not necessarily large),
then an additional term in the RHS of~\eqref{putz} should be added:
\begin{equation}
\label{adit} 
-\nabla E_{x_h,x_k} \colon \int_0^t\!\!\int (u_h u_k)(y,s)\,dy\,ds
\end{equation}
(the $\colon$ notation means that the $\sum_{h,k}$ symbol has been omitted).
Such term is bounded by $C|x|^{-4}t^{1/2}$, thus justifying its inclusion
inside the remainder when $|x|>\!\!\!>\sqrt t>\!\!\!>1$.
With this additional term, the condition of the remainder can be simplified
into
$\lim_{\frac{|x|}{\sqrt t}\to\infty} \frac{\widetilde{\mathcal{R}}(x,t)}{t |x|^{-4}}=0$.
\end{remark}

\begin{remark} [on Corollary~\ref{coro1}]
The restriction $r+\frac3p<a$ in Part~(a) of Corollary~\ref{coro1}
is natural beacause
the decay assumption on $u_0$ guarantees that $e^{t\Delta}u_0\in L^p_r$
exactly for those $r,p$ satisfying such restriction.
However, estimates~\eqref{sdab} remain valid in the whole range $0\le r+\frac3p<3$
(and under the conditions of Part~(b) even for $0\le r<\frac3p<4$
if we want to estimate  $\|u(t)-e^{t\Delta}u_0\|_{L^p_r}$ instead of~$\|u(t)\|_{L^p_r}$
in that expression.
\end{remark}


\bibliographystyle{amsplain}

\end{document}